\documentclass[a4paper,leqno,draft]{article}
\usepackage{amssymb,amsmath,amsfonts,amsthm,
mathrsfs}

\newtheorem{theorem}{Theorem}[section]
\newtheorem{corollary}[theorem]{Corollary}

\newtheorem{proposition}[theorem]{Proposition}
\newtheorem{definition}[theorem]{Definition}
\theoremstyle{definition}
\newtheorem{remark}[theorem]{Remark}

\newcommand{\wt}[1]{\widetilde{#1}}

\newcommand{\Cinf}{\ensuremath{\mathcal{C}^\infty}}
\newcommand{\Cinfc}{\ensuremath{\mathcal{C}^\infty_{\text{c}}}}
\newcommand{\D}{\ensuremath{{\cal D}}}
\renewcommand{\S}{\mathscr{S}}
\newcommand{\E}{\ensuremath{{\cal E}}}
\newcommand{\OM}{\ensuremath{{\cal O}_\mathrm{M}}}

\newcommand{\LL}{\mathcal{L}}

\newcommand{\mb}[1]{\ensuremath{\mathbb{#1}}}
\newcommand{\N}{\mb{N}}

\newcommand{\R}{\mb{R}}
\newcommand{\C}{\mb{C}}

\newcommand{\G}{\ensuremath{{\cal G}}}
\newcommand{\Gt}{\ensuremath{{\cal G}_\tau}}

\newcommand{\Gc}{\ensuremath{{\cal G}_\mathrm{c}}}
\newcommand{\Gcinf}{\ensuremath{{\cal G}^\infty_\mathrm{c}}}
\newcommand{\GS}{\G_{{\, }\atop{\hskip-4pt\scriptstyle\S}}\!}
\newcommand{\EM}{\ensuremath{{\cal E}_{M}}}
\newcommand\EcM{\mathcal{E}_{\mathrm{c},M}}
\newcommand\EcMinf{\mathcal{E}^\infty_{\mathrm{c},M}}
\newcommand\Nc{\mathcal{N}_{\mathrm{c}}}
\newcommand{\Et}{\ensuremath{{\cal E}_{\tau}}}
\newcommand{\ES}{\mathcal{E}_{\S}}
\newcommand{\EMinf}{\ensuremath{{\cal E}^\infty_{M}}}
\newcommand{\ESinf}{\mathcal{E}_{\S}^{\infty}}

\newcommand{\Nt}{\ensuremath{{\cal N}_{\tau}}}
\newcommand{\Neg}{\mathcal{N}}
\newcommand{\NS}{\mathcal{N}_{\S}}

\newcommand{\Ginf}{\ensuremath{\G^\infty}}
\newcommand{\Gtinf}{\mathcal{G}^{\infty}_\tau}
\newcommand{\GSinf}{\G^\infty_{{\, }\atop{\hskip-3pt\scriptstyle\S}}}

\newcommand{\lara}[1]{\langle #1 \rangle}
\newcommand{\WF}{\mathrm{WF}}

\newcommand{\singsupp}{\mathrm{sing\, supp}}
\newcommand{\supp}{\mathrm{supp}}
\newcommand{\Char}{\ensuremath{\text{Char}}}
\newcommand{\zs}{\setminus 0}
\newcommand{\CO}[1]{\ensuremath{T^*(#1) \zs}}
\newcommand{\ssc}{\mathrm{sc}}

\newcommand{\Ellsc}{\mathrm{Ell}_\ssc}

\newfont{\bigmath}{cmr12 at 13pt}
\newcommand{\PPsi}{\bigmath{\symbol{9}}}
\newcommand{\Oprop}[1]{{\ }_{\mathrm{pr}}^{\hphantom{m}}\text{\PPsi}_{\ssc}^{ #1}}
\newcommand{\Opropr}[1]{{\ }_{\mathrm{pr}}^{\hphantom{m}}\text{\PPsi}_{\rm{rg}}^{ #1}}

\newfont{\grecomath}{cmmi12 at 15pt}
\newcommand{\Bigmu}{\text{\grecomath{\symbol{22}}}}

\newcommand{\val}{\mathrm{v}} 
\newcommand{\esp}{\mathrm{e}}

\newfont{\bl}{msbm10 scaled \magstep2}

\newcommand{\beq}{\begin{equation}}
\newcommand{\eeq}{\end{equation}}

\newcommand{\notmid}{\mid\kern-0.5em\not\kern0.5em}

\newcommand{\Ga}{\Gamma}

\newcommand{\eps}{\varepsilon}

\newcommand{\Om}{\Omega}

\newcommand{\compl}[1]{{#1}^{\mathrm{c}}}

\newcommand{\Sscu}{\underline{\mathcal{S}}_{\,\ssc}}
\newcommand{\Nu}{\underline{\mathcal{N}}}
\newcommand{\Syscu}{{\wt{\underline{\mathcal{S}}}}_{\,\ssc}}
\newcommand{\Nuinf}{\Nu^{-\infty}}

\newcommand{\Syru}{{\wt{\underline{\mathcal{S}}}}_{\mathrm{rg}}}

\newcommand{\M}{\mathcal{M}}
\newcommand{\mF}{\mathcal{F}}

\newcommand{\mP}{\mathcal{P}}
\newcommand{\mQ}{\mathcal{Q}}

\newcommand{\Lb}{\mathcal{L}_{\rm{b}}}

\newcommand{\dslash}{d\hspace{-0.4em}{ }^-\hspace{-0.2em}}

\begin{document}

\title{{\bf Microlocal analysis in the dual of a Colombeau algebra: generalized wave front sets and noncharacteristic regularity}}

\author{Claudia Garetto \\
Institut f\"ur Technische Mathematik,\\ Geometrie und Bauniformatik,\\
Universit\"at Innsbruck, Austria\\
\texttt{claudia@mat1.uibk.ac.at}\\
}
\date{ }
\maketitle
\begin{abstract} We introduce different notions of wave front set for the functionals in the dual of the Colombeau algebra $\Gc(\Om)$ providing a way to measure the $\G$ and the $\Ginf$- regularity in $\LL(\Gc(\Om),\wt{\C})$. For the smaller family of functionals having a ``basic structure'' we obtain a Fourier transform-characterization for this type of generalized wave front sets and results of noncharacteristic $\G$ and $\Ginf$-regularity. 
\end{abstract}
 
\setcounter{section}{-1}
\section{Introduction}
The past decade has seen the emergence of a differential-algebraic theory of generalized functions of Colombeau type
\cite{Biagioni:90, Colombeau:85, Colombeau:92, GFKS:01, NPS:98, O:92, Rosinger:90}
that answered a wealth of questions on solutions to linear and nonlinear partial differential equations
involving non-smooth coefficients and strongly singular data.
Interesting results were obtained in Lie group invariance of generalized functions \cite{DKP:02, KO:00, O:01, O:04b},
nonlinear hyperbolic equations with generalized function data
\cite{ColMO:90, MVN:98f, MVN:98, O:87, O:92b, O:04, OW:94, OW:99},
distributional metrics in general relativity \cite{KOSV:04, KS:99, KS:99b, KS:02, KS:02b, KSV:05},
propagation of strong singularities in linear hyperbolic equations with discontinuous coefficients
\cite{HdH:01, HdH:01c, LO:91, O:89}, microlocal analysis, pseudodifferential operators
and Fourier integral operators with non-smooth symbols
\cite{Garetto:04, GGO:03, GH:05, GH:05b, GH:03, GH:04, HO:03, HOP:05}.

Some ``key technologies'' for the regularity theory of partial differential equations in the Colombeau context have been developed in \cite{Garetto:04, Garetto:04th, GGO:03, GH:05}. They consist in a complete theory of generalized pseudodifferential operators (including a parametrix-construction for operators with generalized hypoelliptic symbol) \cite{Garetto:04, GGO:03} and the application of those pseudodifferential techniques to the microlocal analysis of generalized functions \cite{GH:05}. Particular attention has been given to the dual of a Colombeau algebra which plays a main role in the kernel theory for generalized pseudodifferential operators \cite{Garetto:04, GGO:03}. It is now natural to extend the pseudodifferential operator's action to the dual and to shift the microlocal investigations from the level of generalized functions to the level of $\wt{\C}$-linear functionals. This will require notions of local and microlocal regularity in the dual of a Colombeau algebra and suitable ways of measuring such kinds of regularity. Microlocal analysis is essential for a full understanding of the generalized pseudodifferential operator's action and propagation of singularities and has to be developed in the dual context since the kernels of such operators are not always Colombeau generalized functions but functionals. 

The aim of this paper is to provide tools of microlocal analysis suited to investigate the dual of a Colombeau algebra. Based on the duality theory developed within the Colombeau frameworks in \cite{Garetto:05a, Garetto:05b} it extends and adapts the microlocal results for Colombeau generalized functions stated in \cite{GH:05}. In the usual Colombeau context a generalized function $u\in\G(\Om)$ is said to be regular if it belongs to the subalgebra $\Ginf(\Om)$. This allows to set up a regularity theory for $\G(\Om)$ which is coherent with the usual concept of regularity for distributions since $\Ginf(\Om)\cap \D'(\Om)=\Cinf(\Om)$. In \cite{DPS:98, NPS:98} a notion of generalized wave front set is defined for $u\in\G(\Om)$ as a $\Ginf$-wave front set. It means that the conic regions of ``microlocal regularity'' we deal with in the cotangent space are regions of $\Ginf$-regularity. Coming now to the dual $\LL(\Gc(\Om),\wt{\C})$, i.e. the space of all continuous and $\wt{\C}$-linear functionals on $\Gc(\Om)$, where $\wt{\C}$ is the ring of complex generalized numbers, by continuous embedding it contains both $\Ginf(\Om)$ and $\G(\Om)$. As a consequence, two levels of regularity concern a functional in $\LL(\Gc(\Om),\wt{\C})$: the regularity with respect to $\G(\Om)$ and the regularity with respect to $\Ginf(\Om)$. In this paper in order to measure such different kinds of regularity of $T\in\LL(\Gc(\Om),\wt{\C})$ we introduce the notions of $\G$-wave front set ($\WF_\G(T)$) and $\Ginf$-wave front set ($\WF_{\Ginf}(T)$).

Inspired by \cite{GH:05} and making use of the theory of pseudodifferential operators with generalized symbols elaborated in \cite{GGO:03, GH:05}, $\WF_\G(T)$ and $\WF_\Ginf(T)$ are defined as intersection of suitable regions of generalized non-ellipticity of those pseudodifferential operators which map $T$ in $\G(\Om)$ and $\Ginf(\Om)$ respectively. Core of the paper is a Fourier transform-characterization of $\WF_\G(T)$ and $\WF_{\Ginf}(T)$ as in \cite[Theorem 8.56]{Folland:95} which consists in the direct investigation of the properties of the Fourier transform of $T$ after multiplication by a suitable cut-off function. For this purpose special spaces of generalized functions with rapidly decreasing behavior on a conic subset of $\R^n$ are introduced and among all the functionals of $\LL(\Gc(\Om),\wt{\C})$ we restrict to consider those elements which have a ``basic structure''. More precisely we assume $T\in\LL(\Gc(\Om),\wt{\C})$ being defined by a net of distributions $(T_\eps)_\eps$ which fulfills a continuity assumption uniform with respect to $\eps$ (Definition \ref{def_repr_induc}) and the equality $Tu=[(T_\eps u_\eps)_\eps]\in\wt{\C}$ for all $u=[(u_\eps)_\eps]\in\Gc(\Om)$. Even though the $\G$-wave front set and the $\Ginf$-wave front set can be defined on any functional of the dual $\LL(\Gc(\Om),\wt{\C})$, the main theorems and propositions presented here are proven to be valid for basic functionals. In addition, all the results of microlocal regularity have a double version: the $\G$-version and the $\Ginf$-version.

We now describe in detail the contents of the sections.

Section \ref{section_basic} provides the needed theoretical background of basic functionals and refer for topological issues to \cite{Garetto:05a, Garetto:05b}. After the first definitions and basic properties, the action of a basic functional on a Colombeau generalized function in two variables is investigated in Subsection \ref{subsection_action}. Together with some results on the composition of a basic functional with an integral operator in Subsection \ref{subsection_comp}, it gives the essential tools for dealing with the convolution of Colombeau generalized functions and functionals in Subsection \ref{subsection_conv}. The algebras of generalized functions here involved are $\Gc(\Om)$, $\G(\Om)$ and $\GS(\R^n)$ while the functionals are elements of the duals $\LL(\Gc(\Om),\wt{\C})$, $\LL(\G(\Om),\wt{\C})$ and $\LL(\GS(\R^n),\wt{\C})$. A regularization of basic functionals is obtained via convolution with a generalized mollifier. Finally in Subsection \ref{subsection_Fourier} we extends the natural notion of Fourier transform on $\GS(\R^n)$ to the dual $\LL(\GS(\R^n),\wt{\C})$ and we study the Fourier transform of a basic functional in $\LL(\G(\Om),\wt{\C})$. 

In the recent Colombeau literature a pseudodifferential operator with generalized symbol is a $\wt{\C}$-linear continuous operator which maps $\Gc(\Om)$ into $\G(\Om)$. Section \ref{section_pseudo} extends the action of such generalized pseudodifferential operator to the duals $\LL(\Gc(\Om),\wt{\C})$ and $\LL(\G(\Om),\wt{\C})$. The extension procedure is obtained via transposition and gives interesting mapping properties concerning the subspaces of basic functionals. A variety of symbols (and amplitudes) is considered: generalized symbols of order $m$ and type $(\rho,\delta)$, regular symbols, slow scale symbols, generalized symbols of order $-\infty$, regular symbols of order $-\infty$ and generalized symbols of refined order (see \cite{GGO:03, GH:05}). A connection is shown to exist between $\G$-regularity, generalized symbols of order $-\infty$ and basic functionals as well as between $\Ginf$-regularity, regular symbols of order $-\infty$ and basic functionals. More precisely we prove that $R$ is an integral operator with kernel in $\G(\Om\times\Om)$ if and only if it is a pseudodifferential operator with generalized amplitude of order $-\infty$ and that $R$ is $\G$-regularizing on the basic functionals of $\LL(\G(\Om),\wt{\C})$, in the sense that $RT\in\G(\Om)$ if $T\in\LL(\G(\Om),\wt{\C})$ is basic. Analogously $R$ is an integral operator with kernel in $\Ginf(\Om\times\Om)$ if and only if it is a pseudodifferential operator with regular amplitude of order $-\infty$ and it is $\Ginf$-regularizing on the basic functionals of $\LL(\G(\Om),\wt{\C})$. A $\G$-pseudolocality property is obtained for properly supported pseudodifferential operators with generalized symbols while a $\Ginf$-pseudolocality property is valid when the symbols are regular. Section \ref{section_pseudo} ends by adapting the result of $\Ginf$-regularity in \cite{GGO:03} for pseudodifferential operators with generalized hypoelliptic symbols to the dual context of basic functionals.

A $\G$-microlocal analysis and a $\Ginf$-microlocal analysis for the dual $\LL(\Gc(\Om),\wt{\C})$ are settled and developed in Section \ref{section_WF}. The additional assumption of basic structure on the functional $T$ is employed in Subsection \ref{subsection_1} in proving that the projections on $\Om$ of $\WF_\G(T)$ and $\WF_{\Ginf}(T)$ coincide with the $\G$-singular support and the $\Ginf$-singular support of $T$ respectively. The Fourier transform-characterizations of $\WF_\G(T)$ and $\WF_{\Ginf}(T)$ are the result of the $\G$ and the $\Ginf$-microlocal investigations of pseudodifferential operators elaborated throughout Subsection \ref{subsection_2} in the dual $\LL(\Gc(\Om),\wt{\C})$. Concerning the notion of slow scale micro-ellipticity here employed this has been already introduced in \cite{GH:05} while the concept of generalized microsupport of a generalized symbol in \cite[Definition 3.1]{GH:05} is transformed into $\G$-microsupport and $\Ginf$-microsupport (Definition \ref{def_micro_supp}).

Section \ref{section_nonch} concludes the paper with a theorem on noncharacteristic $\G$ and $\Ginf$-regularity for pseudodifferential operators with slow scale symbols when they act on basic functionals of $\LL(\Gc(\Om),\wt{\C})$. This is  an extension and adaptation to the dual $\LL(\Gc(\Om),\wt{\C})$ of Theorem 4.1 in \cite{GH:05}.

For the advantage of the reader we recall in the sequel some topological issues discussed in \cite{Garetto:05a, Garetto:05b} and we fix some notations.

\subsection{Notions of topology and duality theory for spaces of Colombeau type}
\label{subsection_0}
A topological investigation into spaces of generalized functions of Colombeau type has been initiated in \cite{Garetto:05a, Garetto:05b, Garetto:04th, Scarpalezos:92, Scarpalezos:98, Scarpalezos:00} setting the foundations of duality theory in the recent work on topological and locally convex topological $\wt{\C}$-modules \cite{Garetto:05a, Garetto:05b, Garetto:04th}. Without presenting the technical details of this theoretical construction, we recall that a suitable adaptation of the classical notion of seminorm, called ultra-pseudo-seminorm \cite[Definition 1.8]{Garetto:05a}, allows to characterize a locally convex $\wt{\C}$-linear topology as a topology determined by a family of ultra-pseudo-seminorms. The most common Colombeau algebras can be introduced as $\wt{\C}$-modules of generalized functions based on a locally convex topological vector space $E$. Such a $\wt{\C}$-module $\G_E$ is the quotient of the set
\beq
\label{defME}
 \M_E := \{(u_\eps)_\eps\in E^{(0,1]}:\, \forall i\in I\,\, \exists N\in\N\quad p_i(u_\eps)=O(\eps^{-N})\, \text{as}\, \eps\to 0\}
\eeq
of $E$-moderate nets with respect to the set 
\beq
\label{defNE}
 \Neg_E := \{(u_\eps)_\eps\in E^{(0,1]}:\, \forall i\in I\,\, \forall q\in\N\quad p_i(u_\eps)=O(\eps^{q})\, \text{as}\, \eps\to 0\},
\eeq  
of $E$-negligible nets, and it is naturally endowed with a locally convex $\wt{\C}$-linear topology usually called sharp topology in \cite{NPS:98, Scarpalezos:92, Scarpalezos:98, Scarpalezos:00}. Given a family of seminorms $\{p_i\}_{i\in I}$ on $E$, the sharp topology on $\G_E$ is determined by the ultra-pseudo-seminorms $\mP_i(u):=\esp^{-\val_{p_i}(u)}$, where $\val_{p_i}$ is the \emph{valuation} 
\[
\val_{p_i}([(u_\eps)_\eps]):=\val_{p_i}((u_\eps)_\eps):=\sup\{b\in\R:\ p_i(u_\eps)=O(\eps^b)\ \text{as $\eps\to 0$}\}
\]
(see \cite[Subsection 3.1]{Garetto:05a} for further explanations). Note that valuations and ultra-pseudo-seminorms are defined on $\M_E$ and extended to the factor space $\G_E$ in a second time. It is clear that the ring $\wt{\C}$ of complex generalized numbers is an example of $\G_E$-space obtained by choosing $E=\C$. The valuation and ultra-pseudo-norm on $\wt{\C}$ obtained as above by means of the absolute value on $\C$ are denoted by $\val_{\wt{\C}}$ and $|\cdot|_\esp$ respectively. 

As proved in \cite[Corollary 1.17]{Garetto:05a} for an arbitrary locally convex topological $\wt{\C}$-module $(\G,\{\mQ_j\}_{j\in J})$, a $\wt{\C}$-linear map $T:\G_E\to\G$ is continuous if and only if for all $j\in J$ there exists a finite subset $I_0\subseteq I$ and a constant $C>0$ such that for all $u\in\G_E$
\[
\mQ_j(Tu) \le C \max_{i\in I_0}\mP_i(u).
\]

\subsubsection*{The Colombeau algebras $\G(\Om)$, $\Gc (\Om)$, $\GS(\R^n)$}

The Colombeau algebra $\G(\Om)$ is the $\wt{\C}$-module of $\G_E$-type given by $E=\E(\Om)$. Equipped with the family of seminorms $p_{K,i}(f)=\sup_{x\in K, |\alpha|\le i}|\partial^\alpha f(x)|$ where $K\Subset\Om$, the space $\E(\Om)$ induces on $\G(\Om)$ a metrizable and complete locally convex $\wt{\C}$-linear topology  which is determined by the ultra-pseudo-seminorms $\mP_{K,i}(u)=\esp^{-\val_{p_{K,i}}(u)}$. For coherence with some well-established notations in Colombeau theory we write $\M_{\E(\Om)}=\EM(\Om)$ and $\Neg_{\E(\Om)}=\Neg(\Om)$.

The Colombeau algebra $\Gc(\Om)$ of generalized functions with compact support is topologized by means of a strict inductive limit procedure. More precisely, setting $\G_K(\Om):=\{u\in\Gc(\Om):\, \supp\, u\subseteq K\}$ for $K\Subset\Om$, $\Gc(\Om)$ is the strict inductive limit of the sequence $(\G_{K_n}(\Om))_{n\in\N}$, where $(K_n)_{n\in\N}$ is an exhausting sequence of compact subsets of $\Om$ such that $K_n\subseteq K_{n+1}$. We recall that the space $\G_K(\Om)$ is endowed with the topology induced by $\G_{\mathcal{D}_{K'}(\Om)}$ where $K'$ is a compact subset containing $K$ in its interior. In detail we consider on $\G_K(\Om)$ the ultra-pseudo-seminorms $\mP_{\G_K(\Om),n}(u)=\esp^{-\val_{K,n}(u)}$. Note that the valuation $\val_{K,n}(u):=\val_{p_{K',n}}(u)$ is independent of the choice of $K'$ when acts on $\G_K(\Om)$. As observed in \cite[Subsection 1.2.2]{Garetto:04th} the Colombeau algebra $\Gc(\Om)$ is isomorphic to the factor space $\EcM(\Om)/\Nc(\Om)$ where $\EcM(\Om)$ and $\Nc(\Om)$ are obtained by intersecting $\EM(\Om)$ and $\Neg(\Om)$ with $\cup_{K\Subset\Om}\D_K(\Om)^{(0,1]}$ respectively.

The Colombeau algebra $\GS(\R^n)$ of generalized functions based on $\S(\R^n)$ is obtained as a $\G_E$-module by choosing $E=\S(\R^n)$. It is a Fr\'echet $\wt{\C}$-module according to the topology of the ultra-pseudo-seminorms $\mP_h(u)=\esp^{-\val_{p_h}(u)}$, where $p_h(f)=\sup_{x\in\R^n,|\alpha|\le h}(1+|x|)^h|\partial^\alpha f(x)|$, $f\in\S(\R^n)$, $h\in\N$. In the course of the paper we will use the notations $\ES(\R^n)$ and $\NS(\R^n)$ for the spaces of nets $\M_{\S(\R^n)}$ and $\Neg_{\S(\R^n)}$ respectively.

\subsubsection*{The regular Colombeau algebras $\Ginf(\Om)$, $\Gcinf(\Om)$, $\GSinf(\R^n)$}
Given a locally convex topological space $(E,\{p_i\}_{i\in I})$ the $\wt{\C}$-module $\Ginf_E$ of regular generalized functions based on $E$ is defined as the quotient $\M^\infty_E/\Neg_E$, being 
\[
\M^\infty_E :=\{(u_\eps)_\eps\in E^{(0,1]}:\, \exists N\in\N\,\, \forall i\in I\quad p_i(u_\eps)=O(\eps^{-N})\, \text{as}\ \eps\to 0\}
\]
the set of $E$-regular nets. The moderateness properties of $\M_E^\infty$ allows to define the valuation 
\[
\val^\infty_E ((u_\eps)_\eps):=\sup\{b\in\R:\, \forall i\in I\qquad p_i(u_\eps)=O(\eps^b)\ \text{as $\eps\to 0$}\}
\]
which extends to $\Ginf_E$ and leads to the ultra-pseudo-norm $\mP^\infty_E(u):=\esp^{-\val_E^\infty(u)}$. This topological model is employed in endowing the Colombeau algebras $\Ginf(\Om)$, $\Gcinf(\Om)$, $\GSinf(\R^n)$ and $\Gtinf(\R^n)$ with a locally convex $\wt{\C}$-linear topology. 

We begin by recalling that $\Ginf(\Om)$ is the subalgebra of all elements $u$ of $\G(\Om)$ having a representative $(u_\eps)_\eps$ belonging to the set
\[
\EM^\infty(\Om):=\{(u_\eps)_\eps\in\E[\Om]:\ \forall K\Subset\Om\, \exists N\in\N\, \forall\alpha\in\N^n\quad \sup_{x\in K}|\partial^\alpha u_\eps(x)|=O(\eps^{-N})\ \text{as $\eps\to 0$}\}.
\]
$\Ginf(\Om)$ can be seen as the intersection $\cap_{K\Subset\Om}\Ginf(K)$, where $\Ginf(K)$ is the space of all $u\in\G(\Om)$ having a representative $(u_\eps)_\eps$ satisfying the condition: $\exists N\in\N$ $\forall\alpha\in\N^n$,\ $\sup_{x\in K}|\partial^\alpha u_\eps(x)|=O(\eps^{-N})$. The ultra-pseudo-seminorms $\mP_{\Ginf(K)}(u):=\esp^{-\val_{\Ginf(K)}}$, where $\val_{\Ginf(K)}:=\sup\{b\in\R:\, \forall\alpha\in\N^n\quad \sup_{x\in K}|\partial^\alpha u_\eps(x)|=O(\eps^b)\}$, equips $\Ginf(\Om)$ with the topological structure of a Fr\'echet $\wt{\C}$-module.

The algebra $\Gcinf(\Om)$ is the intersection of $\Ginf(\Om)$ with $\Gc(\Om)$. On $\Ginf_K(\Om):=\{u\in\Ginf(\Om):\, \supp\, u\subseteq K\Subset\Om\}$ we define the ultra-pseudo-norm $\mP_{\Ginf_K(\Om)}(u)=\esp^{-\val^\infty_K(u)}$ where $\val^\infty_K(u):=\val^\infty_{\mathcal{D}_{K'}(\Om)}(u)$ and $K'$ is any compact set containing $K$ in its interior. At this point, given an exhausting sequence $(K_n)_n$ of compact subsets of $\Om$, the strict inductive limit procedure determines on a complete and separated locally convex $\wt{\C}$-linear topology on $\Gcinf(\Om)=\cup_n \Ginf_{K_n}(\Om)$. Clearly $\Gcinf(\Om)$ is isomorphic to the factor space $\EcMinf(\Om)/\Nc(\Om)$ where $\EcMinf(\Om):=\EMinf(\Om)\cap(\cup_{K\Subset\Om}\D_K(\Om)^{(0,1]})$.

Finally $\GSinf(\R^n)$ is the $\wt{\C}$-module of regular generalized functions based on $E=\S(\R^n)$. In coherence with the notations already in use we set $\M^\infty_{\S(\Om)}=\ESinf(\Om)$ for any open subset $\Om$ of $\R^n$.


\subsubsection*{The Colombeau algebras of tempered generalized functions $\Gt(\R^n)$ and $\Gtinf(\R^n)$}
The Colombeau algebra of tempered generalized functions $\Gt(\R^n)$ is defined as $\Et(\R^n)/\Nt(\R^n)$, where $\Et(\R^n)$ is the space  
\[
\{(u_\eps)_\eps\in\OM(\R^n)^{(0,1]}:\ \forall\alpha\in\N^n\, \exists N\in\N\quad \sup_{x\in\R^n}(1+|x|)^{-N}|\partial^\alpha u_\eps(x)|=O(\eps^{-N})\ \text{as $\eps\to 0$}\}
\]
of $\tau$-moderate nets and $\Nt(\R^n)$ is the space 
\[
\{(u_\eps)_\eps\in\OM(\R^n)^{(0,1]}:\ \forall\alpha\in\N^n\, \exists N\in\N\, \forall q\in\N\quad \sup_{x\in\R^n}(1+|x|)^{-N}|\partial^\alpha u_\eps(x)|=O(\eps^{q})\ \text{as $\eps\to 0$}\}
\]
of $\tau$-negligible nets. The subalgebra $\Gtinf(\R^n)$ of regular and tempered generalized functions is the quotient $\Et^{\infty}(\R^n)/\Nt(\R^n)$, where $\Et^{\infty}(\R^n)$ is the set of all $(u_\eps)_\eps\in\OM(\R^n)^{(0,1]}$ satisfying the following condition:
\[
\exists N\in\N\, \forall\alpha\in\N^n\, \exists M\in\N\qquad \sup_{x\in\R^n}(1+|x|)^{-M}|\partial^\alpha u_\eps(x)|=O(\eps^{-N}).
\]

\subsubsection*{The topological duals $\LL(\Gc(\Om),\wt{\C})$, $\LL(\Gc(\Om),\wt{\C})$, $\LL(\GS(\R^n),\wt{\C})$}
Throughout the paper the topological duals $\LL(\Gc(\Om),\wt{\C})$, $\LL(\Gc(\Om),\wt{\C})$, $\LL(\GS(\R^n),\wt{\C})$ are endowed with the corresponding topologies of uniform convergence on bounded subsets. These topologies, denoted by $\beta_b(\LL(\Gc(\Om),\wt{\C}),\Gc(\Om))$, $\beta_b(\LL(\G(\Om),\wt{\C}),\G(\Om))$ and $\beta_b(\LL(\GS(\R^n),\wt{\C}),\GS(\R^n))$ respectively, are determined by the ultra-pseudo-seminorms $\mP_B(T):=\sup_{u\in B}|Tu|_\esp$ with $B$ varying in the family of all bounded subsets of $\Gc(\Om)$, $\G(\Om)$ and $\GS(\R^n)$ respectively. As in the classical functional analysis a subset $B$ of a locally convex $\wt{\C}$-module $(\G,\{\mP_i\}_{i\in I})$ is bounded if and only if every ultra-pseudo-seminorm $\mP_i$ is bounded on $B$, i.e. $\sup_{u\in B}\mP_i(u)<\infty$. With respect to the topologies collected in this subsection and the topology on $\Gt(\R^n)$ introduced in \cite[Example 3.9]{Garetto:05a} we have that the following chain of inclusions 
\[
\Ginf(\Om)\,\subseteq\, \G(\Om)\, \subseteq\, \LL(\Gc(\Om),\wt{\C}),
\]
\[
\Gcinf(\Om)\, \subseteq\, \Gc(\Om)\, \subseteq\, \LL(\G(\Om),\wt{\C}),
\]
\[
\GSinf(\R^n)\, \subseteq\, \GS(\R^n)\, \subseteq\, \Gt(\R^n)\, \subseteq\, \LL(\GS(\R^n),\wt{\C})
\]
are continuous \cite[Theorems 3.1, 3.8]{Garetto:05b}. Moreover $\Om\to\LL(\Gc(\Om),\wt{\C})$ is a sheaf and the dual $\LL(\G(\Om),\wt{\C})$ can be identified with the set of functionals in $\LL(\Gc(\Om),\wt{\C})$ having compact support \cite[Theorem 1.2]{Garetto:05b}.

\section{Duality theory in the Colombeau context: basic maps and functionals}
\label{section_basic} 
This section is devoted to maps and functionals defined on $\wt{\C}$-modules of Co\-lom\-beau type. Before considering topics more related to the duals of the Colombeau algebras $\Gc(\Om)$, $\G(\Om)$ and $\GS(\R^n)$ in Subsections \ref{subsection_action}, \ref{subsection_comp}, \ref{subsection_conv} and \ref{subsection_Fourier} we focus our attention on the set $\LL(\G_E,\G_F)$ of all $\wt{\C}$-linear and continuous maps from $\G_E$ to $\G_F$. Among all the elements of $\LL(\G_E,\G_F)$ we study those elements whose action has a ``basic structure'' at the level of representatives.
\begin{definition}
\label{def_repr_map}
Let $(E,\{p_i\}_{i\in I})$ and $(F,\{q_j\}_{j\in J})$ be locally convex topological vector spaces. We say that $T\in\LL(\G_E,\G_F)$ is \emph{basic} if there exists a net $(T_\eps)_\eps$ of continuous linear maps from $E$ to $F$ fulfilling the continuity-property  
\beq
\label{cond_mod_repr}
\forall j\in J\, \exists I_0\subseteq I\, \text{finite}\ \exists N\in\N\, \exists\eta\in(0,1]\, \forall u\in E\, \forall\eps\in(0,\eta]\qquad
q_j(T_\eps u)\le \eps^{-N}\max_{i\in I_0}p_i(u),
\eeq
such that $Tu=[(T_\eps(u_\eps))_\eps]$ for all $u\in\G_E$.
\end{definition}
Note that the equality $Tu=[(T_\eps(u_\eps))_\eps]$ holds for all the representatives of $(u_\eps)_\eps$ since \eqref{cond_mod_repr} entails $(T_\eps v_\eps)_\eps\in\M_F$ if $(v_\eps)_\eps\in\M_E$ and $(T_\eps v_\eps)_\eps\in\Neg_F$ if $(v_\eps)_\eps\in\Neg_E$.
\begin{remark}
\label{rem_representative}
\leavevmode
\begin{trivlist}
\item[(i)] If the net $(T'_\eps)_\eps$ satisfies the condition 
\beq
\label{cond_neg_repr}
\forall j\in J\, \exists I_0\subseteq I\, \text{finite}\ \forall q\in\N\, \exists\eta\in(0,1]\, \forall u\in E\, \forall\eps\in(0,\eta]\qquad
q_j(T'_\eps u)\le \eps^{q}\max_{i\in I_0}p_i(u),
\eeq
then $(T_\eps+T'_\eps)_\eps$ defines the map $T$, in the sense that for all $u\in\G_E$,
\[ 
Tu=[(T_\eps u_\eps)_\eps]=[((T_\eps+T'_\eps)(u_\eps))_\eps]\, \text{in $\G_F$}.
\]
Inspired by the established language of moderateness and negligibility in Co\-lom\-be\-au theory we define the space 
$$\M(E,F):=\{(T_\eps)_\eps\in\LL(E,F)^{(0,1]}:\, (T_\eps)_\eps\, \text{satisfies\ \eqref{cond_mod_repr}}\}$$ of \emph{moderate nets} and the space $$\Neg(E,F):=\{(T_\eps)_\eps\in\LL(E,F)^{(0,1]}:\, (T_\eps)_\eps\, \text{satisfies\  \eqref{cond_neg_repr}}\}$$ of \emph{negligible nets}. By the previous considerations it follows that the classes of $\M(E,F)/\Neg(E,F)$ generate maps in $\LL(\G_E,\G_F)$ which are basic. One easily proves that if $E$ is a normed space with dim$E<\infty$ then the space of all basic maps in $\LL(\G_E,\G_F)$ can be identified with the quotient $\M(E,F)/\Neg(E,F)$. Moreover by Proposition 3.22 in \cite{Garetto:05a} it follows that for any normed space $E$ the ultra-pseudo-normed $\wt{\C}$-module $\G_{E'}$ is isomorphic to the set of all basic functionals in $\LL(\G_E,\wt{\C})$. \item[(ii)] Any continuous linear map $t:E\to F$ produces a natural example of basic element of $\LL(\G_E,\G_F)$. Indeed, as observed in \cite[Remark 3.14]{Garetto:05a}, it is sufficient to take the constant net $(t)_\eps$ and the corresponding map $T:\G_E\to\G_F:u\to[(tu_\eps)_\eps]$. 
\item[(iii)] A certain regularity of the basic operator $T\in\LL(\G_E,\G_F)$ can be already viewed at the level of the net $(T_\eps)_\eps$. Indeed, if we assume that $(T_\eps)_\eps$ belongs to the subset $\M^\infty(E,F)$ of $\M(E,F)$ obtained by substituting the string 
\[
\forall j\in J\ \exists I_0\subseteq I\, \text{finite}\ \exists N\in\N
\]
with
\[
\exists N\in\N\ \forall j\in J\ \exists I_0\subseteq I\, \text{finite}
\]
in \eqref{cond_mod_repr}, we have that $T$ maps $\Ginf_E$ into $\Ginf_F$.
\end{trivlist}
\end{remark}
\begin{definition}
\label{def_repr_induc}
Let $E={\rm{span}}(\cup_{\gamma\in\Gamma}\iota_\gamma(E_\gamma))$, $\iota_\gamma:E_\gamma\to E$ be the inductive limit of the locally convex topological vector spaces $(E_\gamma,\{p_{i,\gamma}\}_{i\in I_\gamma})_{\gamma\in\Gamma}$ and $F$ be a locally convex topological vector space. Let $\G=\wt{\C}-{\rm{span}}(\cup_{\gamma\in\Gamma}\iota_\gamma(\G_{E_\gamma}))\subseteq\G_E$ be the inductive limit of the locally convex topological $\wt{\C}$-modules $(\G_{E_\gamma})_{\gamma\in\Gamma}$. We say that $T\in\LL(\G,\G_F)$ is \emph{basic} if there exists a net $(T_\eps)_\eps\in\LL(E,F)^{(0,1]}$ fulfilling the continuity-property  
\begin{multline}
\label{cond_mod_repr_induc}
\forall\gamma\in\Gamma\, \forall j\in J\, \exists I_{0,\gamma}\subseteq I_\gamma\, \text{finite}\ \exists N\in\N\, \exists\eta\in(0,1]\, \forall u\in E_\gamma\, \forall\eps\in(0,\eta]\\
q_j(T_\eps\iota_\gamma(u))\le \eps^{-N}\max_{i\in I_{0,\gamma}}p_{i,\gamma}(u),
\end{multline}
such that $Tu=[(T_\eps(u_\eps))_\eps]$ for all $u\in\G$.
\end{definition}
It is clear that $(T_\eps)_\eps\in\LL(E,F)^{(0,1]}$ defines a basic map $T\in\LL(\G,\G_F)$ if and only if $(T_\eps\circ\iota_\gamma)_\eps$ defines a basic map $T_\gamma\in\LL(\G_{E_\gamma},\G_F)$ such that $T\circ\iota_\gamma=T_\gamma$ for all  $\gamma\in\Gamma$. We recall that nets $(T_\eps)_\eps$ which define basic maps as in Definitions \ref{def_repr_map} and \ref{def_repr_induc} where already considered in \cite{Delcroix:05,DelSca:00} with slightly more general notions of moderateness and different choices of notations and language.

Particular choices of $E$ and $F$ in the lines above yield the following statements:
\begin{itemize}
\item[(i)] a functional $T\in\LL(\G(\Om),\wt{\C})$ is basic if it is of the form $Tu=[(T_\eps u_\eps)_\eps]$, where $(T_\eps)_\eps$ is a net of distributions in $\E'(\Om)$ satisfying the following condition:
\[
\exists K\Subset\Om\, \exists j\in\N\, \exists N\in\N\, \exists\eta\in(0,1]\, \forall u\in\Cinf(\Om)\, \forall\eps\in(0,\eta]\qquad
|T_\eps(u)|\le \eps^{-N}\sup_{x\in K,|\alpha|\le j}|\partial^\alpha u(x)|.
\]
\item[(ii)] A functional $T\in\LL(\Gc(\Om),\wt{\C})$ is basic if it is of the form $Tu=[(T_\eps u_\eps)_\eps]$, where $(T_\eps)_\eps$ is a net of distributions in $\D'(\Om)$ satisfying the following condition: 
\[
\forall K\Subset\Om\, \exists j\in\N\, \exists N\in\N\, \exists\eta\in(0,1]\, \forall u\in\D_K(\Om)\, \forall\eps\in(0,\eta]\qquad
|T_\eps(u)|\le \eps^{-N}\sup_{x\in K,|\alpha|\le j}|\partial^\alpha u(x)|.
\]
\end{itemize}
Note that in analogy with distribution theory there exists a natural multiplication between functionals in $\LL(\Gc(\Om),\wt{\C})$ and generalized functions in $\G(\Om)$ given by 
\[
uT(v)=T(uv),\qquad\qquad\quad v\in\Gc(\Om).
\]
It provides a $\wt{\C}$-linear operator from $\LL(\Gc(\Om),\wt{\C})$ to $\LL(\Gc(\Om),\wt{\C})$ which maps basic functionals into basic functionals. Moreover, if $u\in\Gc(\Om)$ then $uT\in\LL(\G(\Om),\wt{\C})\subseteq\LL(\G(\R^n),\wt{\C})$.

\subsection{Action of basic functionals on generalized functions in two variables}
\label{subsection_action}
In this subsection we study the action of a basic functional $T$ belonging to the duals $\LL(\Gc(\Om),\wt{\C})$, $\LL(\G(\Om),\wt{\C})$ or $\LL(\GS(\R^n),\wt{\C})$ on a generalized function $u(x,y)$ in two variables. Throughout the paper $\pi_1:\Om'\times\Om\to\Om'$ and $\pi_2:\Om'\times\Om\to\Om$ are the projections of $\Om'\times\Om$ on $\Om'$ and $\Om$ respectively. We recall that $V$ is a proper subset of $\Om'\times\Om$ if for all $K'\Subset\Om'$ and $K\Subset\Om$ we have $\pi_2(V\cap\pi_1^{-1}(K'))\Subset\Om$ and $\pi_1(V\cap\pi_2^{-1}(K))\Subset\Om'$. 
\begin{proposition}
\label{prop_action_func}
Let $\Om'$ be an open subset of $\R^{n'}$ and $T$ be a basic functional of $\LL(\Gc(\Om),\wt{\C})$.
\begin{itemize}
\item[(i)] If $u\in\Gc(\Om'\times\Om)$ then $T(u(x,\cdot)):=[(T_\eps(u_\eps(x,\cdot)))_\eps]$ is a well-defined element of $\Gc(\Om')$;
\item[(ii)] if $u\in\Gcinf(\Om'\times\Om)$ then $T(u(x,\cdot))\in\Gcinf(\Om')$;
\item[(iii)] if $u\in\G(\Om'\times\Om)$ and $\supp\, u$ is a proper subset of $\Om'\times\Om$ then $T(u(x,\cdot))$ defines a generalized function in $\G(\Om')$;
\item[(iv)] $\G$ can be replaced by $\Ginf$ in $(iii)$.
\end{itemize}
Let $T$ be a basic functional of $\LL(\G(\Om),\wt{\C})$.
\begin{itemize}
\item[(v)] If $u\in\G(\Om'\times\Om)$ then $T(u(x,\cdot))\in\G(\Om')$;
\item[(vi)] $\G$ can be replaced by $\Ginf$ in $(v)$;
\item[(vii)] if $u\in\G(\Om'\times\Om)$ and $\supp\, u$ is a proper subset of $\Om'\times\Om$ then $T(u(x,\cdot))\in\Gc(\Om')$;
\item[(viii)] $\G$ can be replaced by $\Ginf$ in $(vii)$.
\end{itemize}
Finally, let $T$ be a basic functional of $\LL(\GS(\R^n),\wt{\C})$.
\begin{itemize}
\item[(ix)] If $u\in\Gt(\R^{2n})$ has a representative $(u_\eps)_\eps$ satisfying the condition
\begin{multline}
\label{cond_A}
\forall\alpha\in\N^n\, \forall s\in\N\, \exists N\in\N\\ 
\quad\sup_{x\in\R^n}(1+|x|)^{-N}\sup_{y\in\R^n,|\beta|\le s}(1+|y|)^s|\partial^\alpha_x\partial^\beta_y u_\eps(x,y)|=O(\eps^{-N})\ \text{as $\eps\to 0$}
\end{multline}
then $T(u(x,\cdot))$ is a well-defined element of $\Gt(\R^n)$;
\item[(x)] if $(u_\eps)_\eps$ fulfills the property
\begin{multline}
\label{cond_B}
\exists M\in\N\, \forall\alpha\in\N^n\, \forall s\in\N\, \exists N\in\N\\
\quad\sup_{x\in\R^n}(1+|x|)^{-N}\sup_{y\in\R^n,|\beta|\le s}(1+|y|)^s|\partial^\alpha_x\partial^\beta_y u_\eps(x,y)|=O(\eps^{-M})\ \text{as $\eps\to 0$}
\end{multline}
then $T(u(x,\cdot))\in\Gtinf(\R^n)$.
\end{itemize}
\end{proposition}
\begin{proof}
$(i)-(ii)$ Let $u\in\Gc(\Om'\times\Om)$ and $(u_\eps)_\eps$ be a representative of $u$ such that $\supp\, u_\eps\subseteq K_1\times K_2$, $K_1\Subset\Om'$, $K_2\Subset\Om$ for all $\eps\in(0,1]$. By definition of basic functional there exists a net $(T_\eps)_\eps\in\D'(\Om)^{(0,1]}$, $N\in\N$, $j\in\N$ and $\eta\in(0,1]$ such that
\beq
\label{est_T}
|T_\eps(u_\eps(x,\cdot))|\le \eps^{-N}\sup_{y\in K_2, |\beta|\le j}|\partial^\beta_y u_\eps(x,y)|
\eeq
for all $x\in\Om'$ and for all $\eps\in(0,\eta]$. From \eqref{est_T} it follows immediately that $(u_\eps)_\eps\in\EcM(\Om'\times\Om)$ implies $(T_\eps(u_\eps(x,\cdot)))_\eps\in\EcM(\Om')$, $(u_\eps)_\eps\in\Nc(\Om'\times\Om)$ implies $(T_\eps(u_\eps(x,\cdot)))_\eps\in\Nc(\Om')$ and $(u_\eps)_\eps\in\EcMinf(\Om'\times\Om)$ implies $(T_\eps(u_\eps(x,\cdot)))_\eps\in\EcMinf(\Om')$. To complete the proof that $T(u(x,\cdot))$ is a well-defined generalized function we still have to prove that it does not depend on the choice of the net $(T_\eps)_\eps$ which determines $T$. Let $(T'_\eps)_\eps\in\D'(\Om)^{(0,1]}$ be another net defining $T$ and $\wt{x}$ a generalized point of $\wt{\Om'}_{\rm{c}}$. Since $u(\wt{x},\cdot):=[(u_\eps(x_\eps,\cdot))_\eps]$ belongs to $\Gc(\Om)$ we have that 
\[
((T_\eps-T'_\eps)(u_\eps(x_\eps,\cdot)))_\eps\in\Neg
\]
i.e., the generalized functions $[(T_\eps(u_\eps(x,\cdot)))_\eps]\in\Gc(\Om')$ and $[(T'_\eps(u_\eps(x,\cdot)))_\eps]\in\Gc(\Om')$ have the same point values. By point value theory this means that $((T_\eps-T'_\eps)(u_\eps(x,\cdot)))_\eps\in\Nc(\Om')$.

$(iii)-(iv)$ Let us now assume that $u\in\G(\Om'\times\Om)$ and that $\supp\, u$ is a proper subset of $\Om'\times\Om$. Let $\chi(x,y)$ be a proper smooth function on $\Om'\times\Om$ identically $1$ in a neighborhood of $\supp\, u$. Clearly we can write $\chi u=u$ in $\G(\Om'\times\Om)$. By the previous reasoning we have that for any $\psi\in\Cinfc(\Om')$ the generalized function $\psi(x)T(u(x,\cdot))=[(\psi(x)T_\eps(\chi(x,\cdot)u_\eps(x,\cdot)))_\eps]$ belongs to $\Gc(\Om')$ if $u\in\G(\Om'\times\Om)$ and to $\Gcinf(\Om')$ if $u\in\Ginf(\Om'\times\Om)$. Finally, let $(\Om'_\lambda)_{\lambda\in\Lambda}$ be a locally finite open covering of $\Om'$ with $\overline{\Om'_\lambda}\Subset\Om'$ and $(\psi_\lambda)_{\lambda\in\Lambda}$ be a family of cut-off functions such that $\psi_\lambda=1$ in a neighborhood of $\pi_2(\supp\, u\, \cap\, \pi_1^{-1}(\overline{\Om'_\lambda}))$. One can easily see that $\psi_\lambda(x)T(u(x,\cdot))|_{\Om'_\lambda}\in\G(\Om'_\lambda)$ determines a coherent family of generalized functions for $\lambda$ varying in $\Lambda$ and therefore, by the sheaf properties of $\G(\Om')$ it defines a generalized function $T(u(x,\cdot))$ in $\G(\Om')$ when $u\in\G(\Om'\times\Om)$. Analogously $T(u(x,\cdot))\in\Ginf(\Om')$ if $u\in\Ginf(\Om'\times\Om)$. We use the notation $T(u(x,\cdot))$ since the definition of this generalized function does not depend on $(\Om'_\lambda)_{\lambda\in\Lambda}$ and $(\psi_\lambda)_{\lambda\in\Lambda}$.

The proof of $(v)$ and $(vi)$ is clear arguing at the level of representatives.

$(vii)-(viii)$ When the support of $u$ is proper we can choose a smooth proper function $\chi(x,y)$ identically $1$ in a neighborhood of $\supp\, u$ and a cut-off function $\psi(y)$ identically $1$ in a neighborhood of $\supp\, T$. Hence, $T(u(x,\cdot))=T(\psi(\cdot)u(x,\cdot))=T(\psi(\cdot)u(x,\cdot)\chi(x,\cdot))$, where $\psi(y)u(x,y)\chi(x,y)\in\Gc(\Om'\times\Om)$. By the first and the second assertion of this proposition we have that $T(u(x,\cdot))$ belongs to $\Gc(\Om')$ or $\Gcinf(\Om')$ if $u$ is an element of $\G(\Om'\times\Om)$ or $\Ginf(\Om'\times\Om)$ respectively.

$(ix)-(x)$ Let us consider $T\in\LL(\GS(\R^n),\wt{\C})$ defined by $(T_\eps)_\eps\in\S'(\R^n)^{(0,1]}$. Recall that $(T_\eps)_\eps$ has the following continuity-property:
\[
\exists j\in\N\, \exists N\in\N\, \exists\eta\in(0,1]\, \forall u\in\S(\R^n)\, \forall\eps\in(0,\eta]\qquad
|T_\eps(u)|\le \eps^{-N}\sup_{y\in\R^n, |\beta|\le j}(1+|y|)^j|\partial^\beta u(y)|.
\]
Hence, if $(u_\eps)_\eps$ satisfies \eqref{cond_A} one gets that for all $x\in\R^n$ and for all $\eps$ small enough
\beq
\label{formula_26}
|\partial^\alpha_x T_\eps(u_\eps(x,\cdot))|=|T_\eps(\partial^\alpha_x u_\eps(x,\cdot))|\le \eps^{-N'}\sup_{y\in\R^n, |\beta|\le j}(1+|y|)^j|\partial^\alpha_x\partial^\beta_y u_\eps(x,y)|\le \eps^{-N'-N}(1+|x|)^{N},
\eeq
where $N$ depends on $\alpha$ and $j$. This means that $(T_\eps(u_\eps(x,\cdot)))_\eps\in\Et(\R^n)$. As already observed in the proof of \cite[Proposition 1.2.25]{Garetto:04th} if $(u_\eps)_\eps$ and $(u'_\eps)_\eps$ are two representatives of $u$ both satisfying \eqref{cond_A} then the difference $v_\eps:=u_\eps-u'_\eps$ fulfills the property
\begin{multline*}
\forall\alpha\in\N^n\, \forall s\in\N\, \exists N\in\N\, \forall q\in\N\\ 
\quad\sup_{x\in\R^n}(1+|x|)^{-N}\sup_{y\in\R^n,|\beta|\le s}(1+|y|)^s|\partial^\alpha_x\partial^\beta_y v_\eps(x,y)|=O(\eps^{q})\ \text{as $\eps\to 0$}.
\end{multline*}
As a consequence $(T_\eps(u_\eps(x,\cdot)-u'_\eps(x,\cdot)))_\eps\in\Nt(\R^n)$. In order to claim that $T(u(x,\cdot))$ is a generalized function in $\Gt(\R^n)$ it remains to show that different nets $(T_\eps)_\eps$ and $(T'_\eps)_\eps$ defining $T$ lead to $((T_\eps-T'_\eps)(u_\eps(x,\cdot)))_\eps\in\Nt(\R^n)$. This is due to the fact that for every $\wt{x}\in\wt{\R^n}$, $u(\wt{x},\cdot):=[(u_\eps(x_\eps,\cdot))_\eps]$ is a generalized function in $\GS(\R^n)$ and then by definition of $T$ the net $((T_\eps-T'_\eps)(u_\eps(x_\eps,\cdot)))_\eps$ is negligible. It follows that the tempered generalized functions $[(T_\eps(u_\eps(x,\cdot)))_\eps]$ and $[(T'_\eps(u_\eps(x,\cdot)))_\eps]$ have the same point values. Thus, $((T_\eps-T'_\eps)(u_\eps(x,\cdot)))_\eps\in\Nt(\R^n)$. Finally, arguing as in \eqref{formula_26} it is clear that $T(u(x,\cdot))$ belongs to $\Gtinf(\R^n)$ when $u$ has a representative satisfying \eqref{cond_B}.
\end{proof}
\begin{remark}
\label{remark_after_prop}
By means of the continuous map
\[
\nu:\GS(\R^n)\to\G(\R^n):(u_\eps)_\eps+\NS(\R^n)\to (u_\eps)_\eps+\Neg(\R^n)
\]
the dual $\LL(\G(\R^n),\wt{\C})$ can be embedded into $\LL(\GS(\R^n),\wt{\C})$ as follows:
\beq
\label{emb_duals}
\LL(\G(\R^n),\wt{\C})\to\LL(\GS(\R^n),\wt{\C}):T\to (u\to T(\nu(u))).
\eeq
Indeed, by composition of continuous maps, $u\to T(\nu(u))$ belongs to the dual $\LL(\GS(\R^n),\wt{\C})$ and, taking a cut-off function $\chi\in\Cinfc(\R^n)$ identically $1$ in a neighborhood of $\supp\, T$, if $u\to T(\nu(u))$ is the null functional in $\LL(\GS(\R^n),\wt{\C})$ we get that 
\[
T(u)=T(\chi u)=T(\nu(\chi u))=0\, \text{in $\wt{\C}$}
\]
for all $u\in\G(\R^n)$. This shows that the map in \eqref{emb_duals} is injective. Obviously all the previous considerations hold for basic functionals.
\end{remark}
Before stating the next proposition we recall that every tempered generalized function can be viewed as an element of $\G(\R^n)$ via the map
\[
\nu_\tau:\Gt(\R^n)\to\G(\R^n):(u_\eps)_\eps+\Nt(\R^n)\to (u_\eps)_\eps+\Neg(\R^n).
\]
\begin{proposition}
Let $T$ be a basic functional of $\LL(\G(\R^n),\wt{\C})$. 
\begin{itemize}
\item[(i)] If $u\in\Gt(\R^{2n})$ then $T((\nu_\tau u)(\xi,\cdot))\in\Gt(\R^n)$;
\item[(ii)] if $u\in\Gtinf(\R^{2n})$ then $T((\nu_\tau u)(\xi,\cdot))\in\Gtinf(\R^n)$;
\item[(iii)] if $u\in\Gt(\R^{2n})$ has a representative $(u_\eps)_\eps$ fulfilling the condition
\beq
\label{cond_repr_GS}
\forall\alpha,\beta,\gamma\in\N^n\, \exists N\in\N\qquad \sup_{y\in\R^n,\xi\in\R^n}(1+|y|)^{-N}|\xi^\beta\partial^\alpha_\xi\partial^\gamma_y u_\eps(\xi,y)|=O(\eps^{-N})
\eeq
then $T((\nu_\tau u)(\xi,\cdot))\in\GS(\R^n)$;
\item[(iv)] if $u\in\Gt(\R^{2n})$ has a representative $(u_\eps)_\eps$ fulfilling the condition
\beq
\label{cond_repr_GSinf}
\exists M\in\N\, \forall\alpha,\beta,\gamma\in\N^n\, \exists N\in\N\qquad \sup_{y\in\R^n,\xi\in\R^n}(1+|y|)^{-N}|\xi^\beta\partial^\alpha_\xi\partial^\gamma_y u_\eps(\xi,y)|=O(\eps^{-M})
\eeq
then $T((\nu_\tau u)(\xi,\cdot))\in\GSinf(\R^n)$.
\end{itemize}
\end{proposition}
\begin{proof}
We begin by observing that the generalized function $T((\nu_\tau u)(\xi,\cdot))$ is defined by the net $(S_\eps(\xi))_\eps:=(T_\eps(u_\eps(\xi,\cdot)))_\eps$, where $(u_\eps)_\eps\in\Et(\R^{2n})$ and $(T_\eps)_\eps$ satisfies the following condition:
\beq
\exists K\Subset\R^n\, \exists j\in\N\, \exists N\in\N\, \exists\eta\in(0,1]\, \forall u\in\Cinf(\R^n)\, \forall\eps\in(0,\eta]\qquad
|T_\eps(u)|\le \eps^{-N}\sup_{y\in K,|\beta|\le j}|\partial^\beta u(y)|.
\eeq
Consequently if $(u_\eps)_\eps\in\Et(\R^{2n})$ then for all $\alpha\in\N^n$ there exists $N'\in\N$ such that for all $\eps$ small enough the estimate
\[
|\partial^\alpha S_\eps(\xi)|=|T_\eps(\partial^\alpha_\xi u_\eps(\xi,y))|\le \eps^{-N}\sup_{y\in K,|\beta|\le j}|\partial^\alpha_\xi\partial^\beta_y u_\eps(\xi,y)|\le c\eps^{-N-N'}(1+|\xi|)^{N'}
\]
holds. This proves that $(S_\eps)_\eps\in\Et(\R^n)$. In an analogous way we obtain that $(S_\eps)_\eps\in\Nt(\R^n)$ when $(u_\eps)_\eps\in\Nt(\R^{2n})$ and that $(S_\eps)_\eps\in\E^\infty_\tau(\R^n)$ when $(u_\eps)_\eps\in\E^\infty_\tau(\R^{2n})$. Note that for all $\xi\in\wt{\R^n}$, $u(\wt{\xi},\cdot):=(u_\eps(\xi_\eps,\cdot))_\eps+\Neg(\R^n)\in\G(\R^n)$. Therefore, for $(T_\eps)_\eps$ and $(T'_\eps)_\eps$ different nets defining $T$ and $(u_\eps)_\eps\in\Et(\R^{2n})$ one has that
\[
(S_\eps(\xi_\eps)-S'_\eps(\xi_\eps)):=(T_\eps(u_\eps(\xi_\eps,\cdot))-T'_\eps(u_\eps(\xi_\eps,\cdot)))_\eps
\]
is negligible. Since $\wt{\xi}$ is arbitrary this implies that $(S_\eps-S'_\eps)_\eps\in\Nt(\R^n)$ and completes the proof of $(i)$ and $(ii)$. Let us assume that $u\in\Gt(\R^{2n})$ has a representative fulfilling \eqref{cond_repr_GS}. Then the corresponding net $(S_\eps)_\eps$, which is already known to belong to $\Et(\R^n)$, satisfies the following estimate: 
\begin{multline*}
\sup_{\xi\in\R^n}|\xi^\beta\partial^\alpha S_\eps(\xi)|=\sup_{\xi\in\R^n}|\xi^\beta T_\eps(\partial^\alpha_\xi u_\eps(\xi,y))|\\
\le \eps^{-N}\sup_{\xi\in\R^n}\sup_{y\in K,|\gamma|\le j}|\xi^\beta\partial^\alpha_\xi\partial^\gamma_y u_\eps(\xi,y)| \le \eps^{-N-N'}\sup_{y\in K}(1+|y|)^{N'},
\end{multline*}
uniformly for small values of $\eps$. This means that $T((\nu_\tau u)(\xi,\cdot))\in\GS(\R^n)$. Finally, when $(u_\eps)_\eps$ satisfies \eqref{cond_repr_GSinf} then 
\[
\sup_{\xi\in\R^n}|\xi^\beta\partial^\alpha S_\eps(\xi)|\le \eps^{-N}\sup_{\xi\in\R^n}\sup_{y\in K,|\gamma|\le j}|\xi^\beta\partial^\alpha_\xi\partial^\gamma_y u_\eps(\xi,y)|\le \eps^{-N-M}\sup_{y\in K}(1+|y|)^{N'},
\]
for some $N,N',M\in\N$. Since $N$ and $M$ do not depend on $\alpha,\beta$ we have that $(S_\eps)_\eps\in\ESinf(\R^n)$ and therefore $T((\nu_\tau u)(\xi,\cdot))\in\GSinf(\R^n)$.
\end{proof}
\begin{remark}
\label{remark_Fourier}
As a straightforward application of the previous proposition we consider the action of a basic functional $T\in\LL(\G(\R^n),\wt{\C})$ on $\esp^{-iy\xi}\in\Gtinf(\R^{2n})$. Omitting the notation $\nu_\tau$ for simplicity, we are allowed to claim that 
\[
T(\esp^{-i\cdot\xi}):=(T_\eps(\esp^{-i\cdot\xi}))_\eps+\Nt(\R^n)
\]
is a generalized function in $\Gtinf(\R^n)$.
\end{remark} 
\subsection{Composition of a basic functional with an integral operator}
\label{subsection_comp}
In the sequel for the advantage of the reader we recall the results on integral operators elaborated in \cite[Proposition 2.14]{GGO:03}. They are needed in stating and proving Proposition \ref{prop_func_int}. 
\begin{proposition}
\label{prop_GGO}
Let us consider the expression 
\beq
\label{int_expression}
\int_{\Om'}k(x,y)u(x)\, dx.
\eeq
\begin{itemize}
\item[(i)] If $k\in\G(\Om'\times\Om)$ then \eqref{int_expression} defines a $\wt{\C}$-linear continuous map
\[
u\to \int_{\Om'}k(x,y)u(x)\, dx
\]
from $\Gc(\Om')$ into $\G(\Om)$;
\item[(ii)] if $k\in\Gc(\Om'\times\Om)$ then \eqref{int_expression} defines a $\wt{\C}$-linear continuous map from $\G(\Om')$ into $\Gc(\Om)$;
\item[(iii)] if $k\in\G(\Om'\times\Om)$ has proper support then the integral operator determined by \eqref{int_expression} maps $\Gc(\Om')$ continuously into $\Gc(\Om)$ and can be uniquely extended to a $\wt{\C}$-linear continuous map from $\G(\Om')$ into $\G(\Om)$.
\end{itemize}
\end{proposition}
\begin{proposition}
\label{prop_func_int}
Let $T$ be a basic functional of $\LL(\G(\Om),\wt{\C})$.
\begin{itemize}
\item[(i)] If $k\in\G(\Om\times\Om)$ and $u\in\Gc(\Om)$ then 
\beq
\label{formula_int}
T\biggl(\int_\Om k(x,y)u(x)\, dx\biggr)=\int_{\Om}T(k(x,\cdot))u(x)\, dx.
\eeq
\item[(ii)] If $\supp\, k$ is proper then \eqref{formula_int} holds for all $u\in\G(\Om)$.
\end{itemize}
Let $T$ be a basic functional of $\LL(\Gc(\Om),\wt{\C})$.
\begin{itemize}
\item[(iii)] If $k\in\Gc(\Om\times\Om)$ then \eqref{formula_int} holds for all $u\in\G(\Om)$;
\item[(iv)] if $k\in\G(\Om\times\Om)$ has proper support then \eqref{formula_int} holds for all $u\in\Gc(\Om)$.
\end{itemize}
\end{proposition}
\begin{proof}
$(i)$ By Proposition \ref{prop_GGO}$(i)$ we know that $\int_\Om k(x,y)u(x)\, dx\in\G(\Om)$ and from Proposition \ref{prop_action_func}$(v)$ we have that $T(k(x,\cdot))\in\G(\Om)$. Therefore it has a meaning the action of $T$ on $\int_\Om k(x,y)u(x)\, dx$ and the integral at the right-hand side of \eqref{formula_int}. The equality is clear since at the level of representatives we can write
\[
T_\eps\biggl(\int_\Om k_\eps(x,y)u_\eps(x)\, dx\biggr)=\int_\Om T_\eps(k_\eps(x,\cdot))u_\eps(x)\, dx.
\]
$(ii)$ If $\supp\, k$ is a proper subset of $\Om\times\Om$ then Proposition \ref{prop_action_func}$(vii)$ says that $T(u(x,\cdot))\in\Gc(\Om)$ and from Proposition \ref{prop_GGO}$(iii)$ it follows that $\int_\Om k(x,y)u(x)\, dx$ defines a generalized function in $\G(\Om)$ when $u\in\G(\Om)$. Let us take a cut-off function $\psi$ identically $1$ in a neighborhood of $\supp\, T$ with $\supp\, \psi\subseteq V\subseteq\overline{V}\Subset\Om$ and a cut-off function $\varphi$ identically $1$ in a neighborhood of $\pi_1(\pi_2^{-1}(\overline{V})\cap\supp\, k)$. By the first assertion we obtain that
\begin{multline*}
T\biggl(\int_\Om k(x,y)u(x)\, dx\biggr)=T\biggl(\psi(y)\int_\Om k(x,y)u(x)\, dx\biggr)\\
=T\biggl(\psi(y)\int_\Om k(x,y)\varphi(x)u(x)\, dx\biggr)=\int_\Om T(k(x,\cdot)\psi(\cdot))\varphi(x)u(x)\, dx\\ =\int_\Om T(k(x,\cdot))u(x)\, dx.
\end{multline*}
$(iii)$ Under the assumptions of $T$ being a basic functional of $\LL(\Gc(\Om),\wt{\C})$, $k\in\Gc(\Om\times\Om)$ and $u\in\G(\Om)$, Proposition \ref{prop_action_func}$(i)$ and Proposition \ref{prop_GGO}$(ii)$ yield that $T(k(x,\cdot))\in\Gc(\Om)$ and $\int_\Om k(x,y)u(x)\, dx\in\Gc(\Om)$ respectively. The equality \eqref{formula_int} is immediate looking at the representatives of the objects involved there.

$(iv)$ Finally, take $k\in\G(\Om\times\Om)$ with proper support and $u\in\Gc(\Om)$. We have that  $T(k(x,\cdot))\in\G(\Om)$ (Proposition \ref{prop_action_func}$(iii)$) and $\int_\Om k(x,y)u(x)\, dx\in\Gc(\Om)$ (Proposition \ref{prop_GGO}$(iii)$). Since $k(x,y)u(x)\in\Gc(\Om\times\Om)$ and $T(k(x,\cdot))u(x)=T(k(x,\cdot)u(x))$ in $\Gc(\Om)$ by the third assertion of this proposition we conclude that the equality \eqref{formula_int} holds.
\end{proof} 
\subsection{Convolution of Colombeau generalized functions and functionals}
\label{subsection_conv}
We proceed by studying the convolution between a Colombeau generalized function and a functional in the dual of the algebras $\Gc(\Om)$, $\G(\Om)$ and $\GS(\R^n)$ or more in general the convolution between two functionals. As in distribution theory this kinds of convolutions are possible under suitable assumptions on the supports of the generalized objects involved. Concerning the functionals we will deal with, a main role is played by the additional hypothesis of ``basic structure''.

We begin by considering the Colombeau generalized function in two variables $u(x-y)$.
One can easily prove that:
\begin{itemize}
\item[$(i)$] if $u\in\G(\R^n)$ then $u(x-y)\in\G(\R^{2n})$;
\item[$(ii)$] if $u\in\Gc(\R^n)$ then $u(x-y)\in\G(\R^{2n})$ and its support is proper;
\item[$(iii)$] $(i)$ and $(ii)$ hold with $\Ginf$ in place of $\G$;
\item[$(iv)$] if $u\in\GS(\R^n)$ then $u(x-y)\in\Gt(\R^{2n})$ and has a representative satisfying condition \eqref{cond_A};
\item[$(v)$] if $u\in\GSinf(\R^n)$ then $u(x-y)\in\Gt^\infty(\R^{2n})$ and has a representative satisfying condition \eqref{cond_B}.
\end{itemize} 
\begin{definition}
\label{def_con_func_functional}
Let $T$ be a basic functional in $\LL(\Gc(\R^n),\wt{\C})$ and $u\in\Gc(\R^n)$. The convolution $u\ast T$ is the generalized function in $\G(\R^n)$ defined by 
\beq
\label{def_conv_1}
u\ast T(x)=T(u(x-\cdot)).
\eeq
\end{definition}
Definition \ref{def_con_func_functional} is the combination of the assertion $(ii)$ above with Proposition \ref{prop_action_func}$(iii)$. Formula \eqref{def_conv_1} allows to define the convolution of $u\in\G(\R^n)$ with a basic functional $T\in\LL(\G(\R^n),\wt{\C})$ (Proposition \ref{prop_action_func}$(v)$) and the convolution of $u\in\GS(\R^n)$ with a basic functional $T\in\LL(\GS(\R^n),\wt{\C})$. In this last case we obtain, by the assertion $(iv)$ above and Proposition \ref{prop_action_func}$(ix)$, that $u\ast T$ is a generalized function in $\Gt(\R^n)$. Analogously, by Proposition  \ref{prop_action_func}$(x)$ it follows that $u\ast T\in\Gt^\infty(\R^n)$ when $u\in\GSinf(\R^n)$.
\begin{proposition}
\label{prop_conv}
If $T$ is a basic functional in $\LL(\G(\R^n),\wt{\C})$ and $u\in\GS(\R^n)$ then $u\ast T\in\GS(\R^n)$.
\end{proposition}
\begin{proof}
Since $\LL(\G(\R^n),\wt{\C})\subseteq\LL(\GS(\R^n),\wt{\C})$ we already know that $u\ast T\in\Gt(\R^n)$. If we prove that this tempered generalized function has a representative in $\ES(\R^n)$ then the proof is complete. By definition of $T$ there exists a net $(T_\eps)_\eps\in\E'(\R^n)^{(0,1]}$, a compact subset $K$ of $\R^n$ and a natural number $j$ such that for some $N\in\N$ for all small enough $\eps$ and for all $f\in\Cinf(\R^n)$,
\beq
\label{assertion}
|T_\eps(f)|\le \eps^{-N}\sup_{y\in K,|\gamma|\le j}|\partial^\gamma f(y)|.
\eeq
Combining \eqref{assertion} with the moderateness properties of $(u_\eps)_\eps\in\ES(\R^n)$ we have that $(T_\eps(u_\eps(x-\cdot)))_\eps$ is a net of functions in $\S(\R^n)$ such that
\begin{multline}
\label{formula_GS}
\sup_{x\in\R^n}|x^\beta\partial^\alpha_x T_\eps(u_\eps(x-\cdot))|= \sup_{x\in\R^n}|T_\eps(x^\beta\partial^\alpha_x u_\eps(x-y))|\\
\le \eps^{-N}\sup_{x\in\R^n,y\in K,|\gamma|\le j}|x^\beta\partial^\alpha_x\partial^\gamma_yu_\eps(x-y)|\le \eps^{-N-N'}\sup_{y\in K}(1+|y|)^{|\beta|},
\end{multline}
where $N'$ depends on $\alpha,\beta\in\N^n$ and $j\in\N$ and the parameter $\eps$ is varying in a sufficiently small interval $(0,\eta]$.
\end{proof}
In the next proposition we collect some continuity results. 
We add a subindex ``b'' in the notation of the duals in order to denote the subspaces of basic functionals. $\LL_{\rm{b}}(\Gc(\R^n),\wt{\C})$, $\LL_{\rm{b}}(\G(\R^n),\wt{\C})$ and $\LL_{\rm{b}}(\GS(\R^n),\wt{\C})$ are equipped with the corresponding topologies of uniform convergence on bounded subsets.

\begin{proposition}
\label{prop_cont_conv}
The $\wt{\C}$-bilinear map
\[
(u,T)\to u\ast T
\]
\begin{itemize}
\item[(i)] from $\Gc(\R^n)\times\Lb(\Gc(\R^n),\wt{\C})$ into $\G(\R^n)$,
\item[(ii)] from $\G(\R^n)\times\Lb(\G(\R^n),\wt{\C})$ into $\G(\R^n)$,
\item[(iii)] from $\GS(\R^n)\times\Lb(\G(\R^n),\wt{\C})$ into $\GS(\R^n)$,
\item[(iv)] from $\Gcinf(\R^n)\times\Lb(\Gc(\R^n),\wt{\C})$ into $\Ginf(\R^n)$,
\item[(v)] from $\Ginf(\R^n)\times\Lb(\G(\R^n),\wt{\C})$ into $\Ginf(\R^n)$,
\item[(vi)] from $\GSinf(\R^n)\times\Lb(\G(\R^n),\wt{\C})$ into $\GSinf(\R^n)$,
\end{itemize}
is separately continuous.
\end{proposition}
\begin{proof}
$(i)$ We want to prove that the map
\beq
\label{map(i)1}
\Gc(\R^n)\to\G(\R^n):u\to u\ast T
\eeq
is continuous for fixed $T\in\LL(\Gc(\R^n),\wt{\C})$ basic and that the map
\beq
\label{map(i)2}
\Lb(\Gc(\R^n),\wt{\C})\to\G(\R^n):T\to u\ast T
\eeq
is continuous for fixed $u\in\Gc(\R^n)$.

We recall that the action of $T$ is given by a net $(T_\eps)_\eps\in\D'(\R^n)^{(0,1]}$ fulfilling the following condition:
\begin{multline}
\label{cond_C}
\forall K\Subset\R^n\, \exists j\in\N\, \exists N\in\N\, \exists\eta\in(0,1]\, \forall u\in\D_K(\R^n)\, \forall\eps\in(0,\eta]\\
|T_\eps(u)|\le \eps^{-N}\sup_{y\in K, |\beta|\le j}|\partial^\beta u(y)|.
\end{multline}
Let us consider the restriction of the map in \eqref{map(i)1} to $\G_K(\R^n)$ and a compact subset $L$ of $\R^n$. Since $\supp\, u(x-\cdot)= x-\supp\, u$, if $u\in\G_K(\R^n)$ then $u(x-\cdot)\in\G_{L-K}(\R^n)$ for all $x\in L$. Under the assumption of $K\subseteq\text{int}(K')\subseteq K'\Subset\R^n$ and $\supp\, u_\eps\subseteq K'$ for all $\eps\in(0,1]$, by \eqref{cond_C} it follows that there exist $j,N\in\N$ such that the estimate
\begin{multline*}
\sup_{x\in L, |\alpha|\le i}|\partial^\alpha_x T_\eps(u_\eps(x-\cdot))|=\sup_{x\in L, |\alpha|\le i}|T_\eps(\partial^\alpha_x u_\eps(x-y))|\\
\le \eps^{-N}\sup_{x\in L, |\alpha|\le i}\sup_{y\in L-K', |\beta|\le j}|\partial^\alpha_x\partial^\beta_y u_\eps(x-y)|\le \eps^{-N}\sup_{z\in K', |\gamma|\le i+j}|\partial^\gamma u_\eps(z)|
\end{multline*}
holds for $\eps$ small. This leads to the continuity of $\G_K(\R^n)\to\G(\R^n):u\to u\ast T$ for all $K\Subset\R^n$. Concerning the map in \eqref{map(i)2} we begin by observing that for $u\in\G_K(\R^n)\subseteq\Gc(\R^n)$ and $L\Subset\R^n$ the set $B_{L,i}:=\{ \partial^\alpha_x u(x-\cdot),\ x\in L,\, |\alpha|\le i\}$ is bounded in $\Gc(\R^n)$, because it is contained in $\G_{L-K}(\R^n)$ and it is bounded there. As a consequence we have the estimate
\beq
\label{cont_T_fixed}
\mP_{L,i}(u\ast T)\le \sup_{v\in B_{L,i}}|T(v)|_\esp,
\eeq
showing the continuity of the map $T\to u\ast T$.
 
$(ii)$ By definition of $T$ basic functional of $\LL(\G(\R^n),\wt{\C})$ and of $(T_\eps)_\eps\in\E'(\R^n)$ one has that 
\[
\exists K\Subset\R^n\, \exists j\in\N\, \exists N\in\N\, \exists\eta\in(0,1]\, \forall u\in\Cinf(\R^n)\, \forall\eps\in(0,\eta]\qquad
|T_\eps(u)|\le \eps^{-N}\sup_{y\in K, |\beta|\le j}|\partial^\beta u(y)|.
\]
Hence, the estimate
\[
\sup_{x\in L, |\alpha|\le i}|T_\eps(\partial^\alpha_x u_\eps(x-\cdot))|\le \eps^{-N}\sup_{x\in L, |\alpha|\le i}\sup_{y\in K, |\beta|\le j}|\partial^\alpha_x\partial^\beta_y u_\eps(x-y)|\le \eps^{-N}\sup_{z\in L-K,|\gamma|\le i+j}|\partial^\gamma u_\eps(z)|,
\]
valid for values of $\eps$ close to $0$, proves that the map $\G(\R^n)\to \G(\R^n):u\to u\ast T$ is continuous. In analogy to $(i)$ we have that for fixed $u\in\G(\R^n)$ the set $B_{L,i}:=\{ \partial^\alpha_x u(x-\cdot),\ x\in L,\, |\alpha|\le i\}$ is bounded in $\G(\R^n)$ and the inequality \eqref{cont_T_fixed} holds for $T\in\Lb(\G(\R^n),\wt{\C})$. Thus, the map $\Lb(\G(\R^n),\wt{\C})\to\G(\R^n):T\to u\ast T$ is continuous.

$(iii)$ Proposition \ref{prop_conv} proves that when $T$ is a basic functional in $\LL(\G(\R^n),\wt{\C})$ and $u\in\GS(\R^n)$ then $u\ast T\in\GS(\R^n)$. More precisely, \eqref{formula_GS} yields the inequality
\[
\sup_{x\in\R^n}|x^\beta\partial^\alpha_x T_\eps(u_\eps(x-\cdot))|\le \eps^{-N}\sup_{z\in\R^n, |\gamma|\le j+|\alpha|}(1+|z|)^{|\beta|}|\partial^\gamma u_\eps(z)|,
\]
where $N,j$ depend only on $T$ and $\eps$ is small enough. This means that for all $h\in\N$
\[
\mP_h(u\ast T)\le \esp^N\mP_{h+j}(u).
\]
In other words the map $\GS(\R^n)\to\GS(\R^n):u\to u\ast T$ is continuous. As in $(ii)$, for fixed $u\in\GS(\R^n)$ the set $B_i:=\{\lara{x}^i\partial^\alpha_x u(x-\cdot),\ |\alpha|\le i, x\in\R^n\}$ is bounded in $\G(\R^n)$. Therefore, the continuity of $\Lb(\G(\R^n),\wt{\C})\to \GS(\R^n):T\to u\ast T$ is due to
\[
\mP_i(u\ast T)\le \sup_{v\in B_i}|T(v)|_\esp.
\]

For the sake of brevity we omit the complete proof of assertions $(iv)$, $(v)$ and $(vi)$. We only remark that in proving the continuity of the map $T\to u\ast T$ in $(iv)$ and $(v)$ respectively we make use given $L\Subset\R^n$ of the subset $B_L:=\{\partial^\alpha_x u(x-\cdot),\ x\in L, \alpha\in\N^n\}$ which is bounded in $\Gc(\R^n)$ when $u\in\Gcinf(\R^n)$ and bounded in $\G(\R^n)$ when $u\in\Ginf(\R^n)$. In proving $(vi)$ we employ the bounded subset $B:=\{\lara{x}^i\partial^\alpha_x u(x-\cdot),\ |\alpha|\le i,\, i\in\N,\, x\in\R^n\}$ of $\G(\R^n)$ where $u\in\GSinf(\R^n)$.

\end{proof}
\begin{definition}
\label{def_conv_func}
Let $S\in\LL(\G(\R^n),\wt{\C})$ and $T$ be a basic functional of the dual $\LL(\Gc(\R^n),\wt{\C})$. The convolution $S\ast T$ is a functional in $\LL(\Gc(\R^n),\wt{\C})$ defined by
\beq
\label{formula_conv_func}
S \ast T(u) = S_x(T_y(u(x+y))).
\eeq
\end{definition}
Definition \ref{def_conv_func} is meaningful since \eqref{formula_conv_func} can be rewritten as 
\[
S((\wt{u}\ast T)^{\wt{\ }}),\qquad 
\]
where $\wt{v}(y):=v(-y)$ is a continuous map from $\Gc(\R^n)$ and $\G(\R^n)$ into themselves respectively. By Proposition \ref{prop_cont_conv}$(i)$ the map $\Gc(\R^n)\to\G(\R^n):u\to(\wt{u}\ast T)^{\wt{\ }}$ is continuous and by composition with the $\wt{\C}$-linear and continuous functional $S$ we conclude that $S\ast T\in\LL(\Gc(\R^n),\wt{\C})$. \begin{proposition}
\label{prop_support}
If $S\in\LL(\G(\R^n),\wt{\C})$ and $T$ is a basic functional of the dual $\LL(\Gc(\R^n),\wt{\C})$ then 
\beq
\label{formula_supports}
\supp(S\ast T)\subseteq \supp\, S + \supp\, T.
\eeq
\end{proposition}
\begin{proof}
Let $A=\supp\, T$ and $B=\supp\, S$. Then $A+B$ is a closed subset of $\R^n$. Let $V=\R^n\setminus(A+B)$ and $u\in\Gc(V)$. Since $\supp(u(x+y))\subseteq\{(x,y):\, x+y\in V\}$ then $S\ast T(u)=0$ and the proof is complete.
\end{proof}
Proposition \ref{prop_support} proves that the convolution of $S\in\LL(\G(\R^n),\wt{\C})$ with  $T\in\LL(\G(\R^n),\wt{\C})$ basic is an element of the dual $\LL(\G(\R^n),\wt{\C})$. It is clear that in all the situations considered so far $S\ast T$ is basic if both $S$ and $T$ are basic.
\begin{proposition}
\label{prop_extension}
The convolution product $\ast$ between functionals extends the convolution product between Colombeau generalized functions and functionals.
\end{proposition}
\begin{proof}
Let $T\in\LL(\Gc(\R^n),\wt{\C})$ be basic and $u\in\Gc(\R^n)\subseteq\LL(\G(\R^n),\wt{\C})$. By \eqref{formula_conv_func} for all $v\in\Gc(\R^n)$ we can write
\[
(u\ast T)(v)=\int_{\R^n}u(x)T(v(x+\cdot))\, dx
\]
Proposition \ref{prop_func_int}$(iv)$ leads to 
\[
\int_{\R^n}u(x)T(v(x+\cdot))\, dx = T\biggl(\int_{\R^n}v(x+\cdot)u(x)\, dx\biggr) = T\biggl(\int_{\R^n}v(z)u(z-\cdot)\, dz\biggr)=\int_{\R^n}v(z)T(u(z-\cdot))\, dz
\]
and shows that Definition \ref{def_conv_func} coincides with Definition \ref{def_con_func_functional} on the couple $(u,T)$. In the same way, making use of assertions $(i)$ and $(ii)$ of Proposition \ref{prop_func_int} one can prove that 
\[
(u\ast T)(v)=\int_{\R^n}v(z)T(u(z-\cdot))\, dz\, \qquad\qquad v\in\Gc(\R^n)
\]
when $u\in\G(\R^n)$ and $T$ is a basic functional of $\LL(\G(\R^n),\wt{\C})$. 
\end{proof}
\begin{remark}
\label{remark_imp}
Combining Proposition \ref{prop_cont_conv} with Propositions \ref{prop_support} and \ref{prop_extension} we obtain that if $u\in\Gc(\R^n)$ and $T\in\LL(\G(\R^n),\wt{\C})$ is basic then $u\ast T\in\Gc(\R^n)$. In particular Proposition \ref{prop_cont_conv}$(v)$ yields that $u\ast T\in\Gcinf(\R^n)$ when $u\in\Gcinf(\R^n)$. We leave to the reader to check that in both these cases the convolution product is a separately continuous $\wt{\C}$-bilinear map. It follows that \eqref{formula_conv_func} applies to $S\in\LL(\Gc(\R^n),\wt{\C})$ and $T$ basic functional in $\LL(\G(\R^n),\wt{\C})$ and defines an element of $\LL(\Gc(\R^n),\wt{\C})$.
\end{remark}
\begin{proposition}
\label{prop_cont_conv_functional}
The $\wt{\C}$-bilinear map
\[
(S,T)\to S\ast T
\]
\begin{itemize}
\item[(i)] from $\LL(\G(\R^n),\wt{\C})\times\Lb(\Gc(\R^n),\wt{\C})$ into $\LL(\Gc(\R^n),\wt{\C})$,
\item[(ii)] from $\LL(\Gc(\R^n),\wt{\C})\times\Lb(\G(\R^n),\wt{\C})$ into $\LL(\Gc(\R^n),\wt{\C})$,
\item[(iii)] from $\LL(\G(\R^n),\wt{\C})\times\Lb(\G(\R^n),\wt{\C})$ into $\LL(\G(\R^n),\wt{\C})$,
\end{itemize}
is separately continuous. 
\end{proposition}
\begin{proof}
We begin by writing the action of the convolution product $S\ast T$ on $u$ as $S((\wt{u}\ast T)^{\wt{\,}})$.

$(i)$ We fix a basic functional $T\in\LL(\Gc(\R^n),\wt{\C})$ and a bounded subset $B\subseteq\Gc(\R^n)$. By Proposition \ref{prop_cont_conv}$(i)$ the map $u\to (\wt{u}\ast T)^{\wt{\, }}$ is continuous from $\Gc(\R^n)$ into $\G(\R^n)$, then the set $B':=\{(\wt{u}\ast T)^{\wt{\, }},\, u\in B\}$ is bounded in $\G(\R^n)$. As a consequence, the equality
\[
\sup_{u\in B}|(S\ast T)(u)|_\esp = \sup_{v\in B'}|S(v)|_\esp
\]
holds and proves the continuity of $\LL(\G(\R^n),\wt{\C})\to\LL(\Gc(\R^n),\wt{\C}):S\to S\ast T$.\\
Let us now fix $S\in\LL(\G(\R^n),\wt{\C})$. Since it is continuous there exist some compact set $K$ and a natural number $m$ such that for all $u\in B$, $B$ bounded subset of $\Gc(\R^n)$, one has
\beq
\label{cont_1}
\sup_{u\in B}|(S\ast T)(u)|_\esp\le c\sup_{u\in B}\mP_{K,m}((\wt{u}\ast T)^{\wt{\, }}).
\eeq
Note that 
\beq
\label{cont_2}
\mP_{K,m}((\wt{u}\ast T)^{\wt{\, }})\le \sup_{v\in B'}|T(v)|_\esp
\eeq
where $B':=\{\partial^\alpha_x u(\wt{x}+\cdot),\, u\in B,\, \wt{x}\in \wt{K},\, |\alpha|\le m\}$ is a bounded subset of $\Gc(\R^n)$. This is due to the fact that, working at the level of representatives, we have  
\[
\val\big(\big(\sup_{x\in K,|\alpha|\le m}|T_\eps(\partial^\alpha u_\eps(x+\cdot))|\big)_\eps\big)\ge\min_{|\alpha|\le m}\val\big(\big(T_\eps(\partial^\alpha u_\eps(x_\eps +\cdot))\big)_\eps\big),
\]
where $(x_\eps)_\eps$ is a net of points of $K$. \eqref{cont_1} combined with \eqref{cont_2} proves that the map $\Lb(\Gc(\R^n),\wt{\C})\to \LL(\Gc(\R^n),\wt{\C}):T\to S\ast T$ is continuous.

$(ii)$ We omit the details of the proof since the arguing is analogous to the one adopted in the first case. We only observe that Remark \ref{remark_imp} is employed in proving the desired continuity.  

$(iii)$ Finally we assume that both the functionals $S,T$ belong to $\LL(\G(\R^n),\wt{\C})$ and that $T$ is basic. For $B$ bounded subset of $\G(\R^n)$, by Proposition \ref{prop_cont_conv}$(ii)$ it follows that $B':=\{(\wt{u}\ast T)^{\wt{\, }},\, u\in B\}\subseteq\G(\R^n)$ is bounded. Thus, the equality $\sup_{u\in B}|(S\ast T)(u)|_\esp = \sup_{v\in B'}|S(v)|_\esp$ shows that the map $S\to S\ast T$ is continuous. Consider now $S$ fixed. For some $K\Subset\R^n$ and $m\in\N$ it holds that 
\[
\sup_{u\in B}|(S\ast T)(u)|_\esp\le c\sup_{u\in B}\mP_{K,m}((\wt{u}\ast T)^{\wt{\, }})\le c\sup_{v\in B'}|T(v)|_\esp,
\]
where $B':=\{\partial^\alpha_x u(\wt{x}+\cdot),\, u\in B,\, \wt{x}\in \wt{K},\, |\alpha|\le m\}\subseteq\G(\R^n)$ is bounded.
\end{proof}
We conclude this subsection with the following regularization of basic functionals.
\begin{theorem}
\label{theorem_reg}
Let $\rho\in\Cinfc(\R^n)$ with $\int_{\R^n}\rho(y) dy=1$. 
\begin{itemize}
\item[(i)] If $T$ is a basic functional of $\LL(\Gc(\R^n),\wt{\C})$ then $[(\rho_{\eps^q})_\eps]\ast T\in\G(\R^n)$ for all $q\in\N$ and $$[(\rho_{\eps^q})_\eps]\ast T\to T$$ in $\LL(\Gc(\R^n),\wt{\C})$ as $q\to \infty$.
\item[(ii)] If $T$ is a basic functional of $\LL(\G(\R^n),\wt{\C})$ then $[(\rho_{\eps^q})_\eps]\ast T\to T$ in $\LL(\G(\R^n),\wt{\C})$ as $q\to \infty$.
\item[(iii)] If $T$ is a basic functional of $\LL(\GS(\R^n),\wt{\C})$ then $[(\rho_{\eps^q})_\eps]\ast T\to T$ in $\LL(\GS(\R^n),\wt{\C})$ as $q\to \infty$.
\end{itemize}
\end{theorem}
\begin{proof}
$(i)$ We begin by proving that when $u\in\Gc(\R^n)$ then $u_q:=[(\rho_{\eps^q})_\eps]\ast u\to u$ in $\Gc(\R^n)$ and that this convergence is uniform on bounded subsets of $\Gc(\R^n)$. If $u\in\Gc(\R^n)$ then $u\in\G_K(\R^n)$ for some $K\Subset\R^n$ and by Proposition \ref{prop_support} $(u_q)_q$ is a sequence of generalized functions in $\G_{K_1}(\R^n)$ with $K_1=K+B_r(0)$, $\supp\, \rho\subseteq B_r(0):=\{y\in\R^n:\, |y|\le r\}$. In particular, at the level of representatives one has that 
\begin{multline*}
\partial^\beta(\rho_{\eps^q}\ast u_\eps -u_\eps)(x)= \int_{\R^n}\rho(z)[\partial^\beta u_\eps(x-\eps^q z)-\partial^\beta u_\eps(x)]\, dz\\
=\int_{\R^n}\rho(z)\sum_{|\alpha|=1}\frac{1}{\alpha!}\partial^{\alpha+\beta}u_\eps(x-\eps^q\theta z)(-\eps^q z)\, dz
\end{multline*}
and for some $K'_1\Subset\R^n$ such that $K_1\subset{\text{int}}(K'_1)\subseteq K'_1$,
\beq
\label{convergence_1}
\sup_{x\in K'_1}|\partial^\beta(\rho_{\eps^q}\ast u_\eps -u_\eps)(x)|\le c\eps^{q-N},
\eeq
where $N$ depends only on $\beta$, $u$ and $\rho$. \eqref{convergence_1} yields that $u_q\to u$ in $\G_{K_1}(\R^n)$. This convergence is uniform on bounded subsets of $\Gc(\R^n)$. Indeed, if $B\subseteq\Gc(\R^n)$ is bounded then it is contained in some $\G_K(\R^n)$ and bounded there. Thus, for some $K'\Subset\R^n$ with $K\subseteq{\text{int}}(K')\subseteq K'$ by the previous computations the inequality 
\[
\sup_{x\in K'_1}|\partial^\beta(\rho_{\eps^q}\ast u_\eps -u_\eps)(x)|\le c\eps^q\sup_{|\gamma|\le|\beta|+1, y\in K'}|\partial^\gamma u_\eps(y)|
\]
holds and leads to $\val_{K_1,|\beta|}(u_q-u)\ge q+\val_{K,|\beta|+1}(u)$. By the assumption of boundedness of $B$ we have that there exists $N\in\N$ such that $\val_{K,|\beta|+1}(u)\ge -N$ for all $u\in B$. Hence, $\sup_{u\in B}\mP_{\G_{K_1}(\R^n),|\beta|}(u_q-u)\ge q-N$ or in other words $$\sup_{u\in B}\mP_{\G_{K_1}(\R^n),|\beta|}(u_q-u)\to 0.$$
Let us now consider $T\in\LL(\Gc(\R^n),\wt{\C})$ basic and $T_q:=[(\rho_{\eps^q})_\eps]\ast T\in\G(\R^n)$. For all $u\in\Gc(\R^n)$ since $(T_\eps\ast \rho_{\eps^q})\ast \wt{u_\eps}(0)-T_\eps\ast\wt{u_\eps}(0)=(T_\eps\ast(\rho_{\eps^q}\ast\wt{u_\eps}-\wt{u_\eps}))(0)$, where $\wt{u}(x)=u(-x)$, we have that
\beq
\label{imp_eq}
(T_q-T)(u)=T((\wt{u}_q-\wt{u})^{\wt{\, }}).
\eeq
By the uniform convergence proved above we know that when $B$ is a bounded subset of $\Gc(\R^n)$ there exists $K_1\Subset\R^n$ such that
\beq
\label{convergence_2}
\forall m\in\N\, \forall\eta\in(0,1]\, \exists\overline{q}\in\N\, \forall q\ge\overline{q}\qquad \sup_{u\in B}\mP_{\G_{K_1}(\R^n),m}((\wt{u}_q-\wt{u})^{\wt{\, }})\le\eta.
\eeq
\eqref{convergence_2} combined with the continuity of $T$ implies that for all $q\ge\overline{q}$ 
\[
\sup_{u\in B}|(T_q-T)(u)|_\esp\le c\sup_{u\in B}\mP_{\G_{K_1}(\R^n),m}((\wt{u}_q-\wt{u})^{\wt{\, }})\le c\,\eta.
\]
This means that $T_q\to T$ in $\LL(\Gc(\R^n),\wt{\C})$ according to $\beta_b(\LL(\Gc(\R^n),\wt{\C}),\Gc(\R^n))$.

$(ii)$ Let now $T$ be a basic functional of $\LL(\G(\R^n),\wt{\C})$. When $u\in\G(\R^n)$ the sequence $u_q:=[(\rho_{\eps^q})_\eps]\ast u$ converges to $u$ uniformly on bounded subsets of $\G(\R^n)$. Indeed, from the same computations of case $(i)$ we have that
\beq
\label{convergence_3}
\sup_{u\in B}\mP_{K,|\beta|}(u_q-u)\le \esp^{-q}\sup_{u\in B}\mP_{K_1,|\beta|+1}(u),
\eeq
where $K_1=K+B_r(0)$. Since $T$ is continuous from $\G(\R^n)$ to $\wt{\C}$ by \eqref{convergence_3} and the equality \eqref{imp_eq} we conclude that $T_q\to T$ with respect to the topology $\beta_b(\LL(\G(\R^n),\wt{\C}),\G(\R^n))$.

We omit the proof of the third assertion since it consists in showing that $u_q$ tends to $u$ uniformly on any bounded subset of $\GS(\R^n)$.
\end{proof}
For the sake of completeness note that the previous regularization for basic functionals in $\LL(\GS(\R^n),\wt{\C})$ is valid for $\rho\in\S(\R^n)$ with $\int\rho(x)\,dx=1$. Such a convergence result was already stated in \cite[Proposition 3.12]{Garetto:05b} for a special family of basic functionals: the generalized delta functionals $\delta_{\wt{x}}$ with $\wt{x}\in\wt{\R^n}$.

\subsection{Fourier transform in the dual $\LL(\GS(\R^n),\wt{\C})$}
\label{subsection_Fourier}
We conclude this section by introducing a natural notion of Fourier transform in the dual of the Colombeau algebra $\GS(\R^n)$. We recall that the Fourier transform $\mF u$ and the inverse Fourier transform $\mF^\ast u$ of a generalized function $u\in\GS(\R^n)$ are defined by the corresponding transformations at the level of representatives. $\mF$ and $\mF^\ast$ are continuous isomorphisms from $\GS(\R^n)$ onto $\GS(\R^n)$. Note that the continuity is the consequence of $\mF$ and $\mF^\ast$ being basic maps on $\GS(\R^n)$. Finally, as explained in detail in \cite[Subsection 1.2.6]{Garetto:04th} all the properties which hold for the transformations on $\S(\R^n)$ can be stated on $\GS(\R^n)$.
\begin{definition}
\label{def_Fourier_dual}
Let $T\in\LL(\GS(\R^n),\wt{\C})$. We define the \emph{Fourier transform} of $T$ as the functional $\mF T$ on $\GS(\R^n)$ given by the formula
\beq
\label{formula_Fourier}
\mF(T)(u)=T(\mF u).
\eeq
\end{definition} 
By the continuity of $\mF:\GS(\R^n)\to\GS(\R^n)$ it is clear that $\mF T\in\LL(\GS(\R^n),\wt{\C})$. Moreover the embedding of $\GS(\R^n)$ into $\LL(\GS(\R^n),\wt{\C})$ proved in \cite[Theorem 3.8]{Garetto:05b} shows that $\mF$ is an extension of the Fourier transform on $\GS(\R^n)$ to $\LL(\GS(\R^n),\wt{\C})$. This motivates the choice of the same notation on $\GS(\R^n)$ and its dual. Obviously the inverse Fourier transform is defined by replacing $\mF$ with $\mF^\ast$ in \eqref{formula_Fourier}. $\mF$ and $\mF^\ast$ are continuous isomorphisms on $\LL(\GS(\R^n),\wt{\C})$. In fact if $B$ is a bounded subset of $\GS(\R^n)$ then $\mF(B):=\{\mF(u),\, u\in B\}$ is bounded itself and the equality
\[
\sup_{u\in B}|\mF T(u)|_\esp=\sup_{u\in B}|T(\mF u)|_\esp=\sup_{v\in\mF(B)}|Tv|_\esp
\]
holds. Analogously this kind of arguing is valid for $\mF^\ast$. We leave to the reader to verify that basic functionals are mapped into basic functionals by $\mF$ and $\mF^\ast$.

In Section \ref{section_WF} we will often deal with the Fourier transform on a basic functional of $\LL(\G(\R^n),\wt{\C})$. By combining some of the results presented in Subsection \ref{subsection_action} we arrive at the following conclusion.
\begin{proposition}
\label{prop_Fourier_comp}
The Fourier transform of a basic functional in $\LL(\G(\R^n),\wt{\C})$ is the tempered generalized function in $\Gtinf(\R^n)$ given by $T(\esp^{-i\cdot\xi}):=(T_\eps(\esp^{i\cdot\xi}))_\eps+\Nt(\R^n)$.
\end{proposition}
\begin{proof}
The functional $T$ can be expressed by a net $(T_\eps)_\eps$ of distributions in $\E'(\R^n)$. Hence, denoting the Fourier transform on $\S'(\R^n)$ by $\widehat{\,}$ we have that $\mF T(u)=[(\widehat{T_\eps}(u_\eps))_\eps]$ and by classical arguments the equality
\[
\widehat{T_\eps}(u_\eps)=\int_{\R^n}T_\eps(\esp^{-i\cdot\xi})u_\eps(\xi)\, d\xi
\]
 is valid for all values of $\eps$. By Remark \ref{remark_Fourier} it follows that for all $u\in\GS(\R^n)$, $$\mF T(u)=\int_{\R^n}T(\esp^{-i\cdot\xi})u(\xi)\, d\xi,$$ where $T(\esp^{-i\cdot\xi}):=(T_\eps(\esp^{i\cdot\xi}))_\eps+\Nt(\R^n)$.
\end{proof}

\section{Generalized pseudodifferential operators acting on the duals $\LL(\Gc(\Om),\wt{\C})$ and $\LL(\G(\Om),\wt{\C})$}
\label{section_pseudo}
In the Colombeau literature a systematic approach to the theory of generalized pseudodifferential operators is given for the first time in \cite{GGO:03}. Based on a notion of generalized symbol as equivalence class it develops a full local calculus for the corresponding pseudodifferential operators acting on the Colombeau algebras $\Gc(\Om)$ and $\G(\Om)$. Results of $\Ginf$-regularity are obtained by means of a parametrix construction for a certain family of operators whose generalized symbols satisfy suitable hypoellipticity assumptions. Concerning this issue a main role is played by different scales in $\eps$ at the level of representatives and by the concept of \emph{slow scale net}. We say that $(\omega_\eps)_\eps\in\C^{(0,1]}$ is a slow scale net if for all $p\ge 0$ there exists $c_p>0$ such that $|\omega_\eps|^p\le c_p\eps^{-1}$ for all $\eps\in(0,1]$. Sometimes the additional assumption of $\inf_\eps\omega_\eps\ge c>0$ is required on $\omega_\eps\in\R^{(0,1]}$ for technical reasons. In this case $\omega_\eps\in\R^{(0,1]}$ is said to be a \emph{strongly positive slow scale net}. In the sequel $\Pi_\ssc$ denotes the set of all strongly positive slow scale nets. Finally we recall that $S^m_{\rho,\delta}(\Om\times\R^p)$ is the usual set of H\"ormander symbols of order $m$ and type $(\rho,\delta)$, with $\rho\in(0,1]$, $\delta\in[0,1)$ and $\Om$ open subset of $\R^n$. $S^m_{\rho,\delta}(\Om\times\R^p)$ is a Fr\'{e}chet space endowed with the seminorms 
\[ 
 |a|^{(m)}_{\rho,\delta,K,\alpha,\beta} = \sup_{x\in K,\xi\in\R^p}\langle\xi\rangle^{-m+\rho|\alpha|-\delta|\beta|}|\partial^\alpha_\xi\partial^\beta_x a(x,\xi)|,
\]
where $K$ ranges over the compact subsets of $\Om$.

In this paper we make use of the following sets of symbols:
\begin{itemize}
\item[-] \emph{generalized symbols}:\qquad $\wt{\mathcal{S}}^m_{\rho,\delta}(\Om\times\R^p):=\G_{S^m_{\rho,\delta}(\Om\times\R^p)}$,
\item[-] \emph{regular symbols}:\qquad $\wt{\mathcal{S}}^m_{\rho,\delta,\rm{rg}}(\Om\times\R^p):=\displaystyle\cap_{K\Subset\Om}	\,\wt{\mathcal{S}}^m_{\rho,\delta,\rm{rg}}(K\times\R^p)$,
\item[-] \emph{slow scale symbols}:\qquad $\wt{\mathcal{S}}^m_{\rho,\delta,\ssc}(\Om\times\R^p):=\mathcal{S}_{\rho,\delta,\ssc}(\Om\times\R^p)/\Neg_{\rho,\delta}(\Om\times\R^p)$,
\item[-] \emph{generalized symbols of order $-\infty$}:\qquad $
\wt{\mathcal{S}}^{-\infty}(\Om\times\R^p):=\G_{S^{-\infty}(\Om\times\R^p)}$,
\item[-] \emph{regular symbols of order $-\infty$}:\qquad $\wt{\mathcal{S}}^{-\infty}_{\rm{rg}}(\Om\times\R^p):=\cap_{K\Subset\Om}\,\wt{\mathcal{S}}^{-\infty}_{\rm{rg}}(K\times\R^p)$,
\end{itemize}
where, for $K$ compact subset of $\Om$, $\wt{\mathcal{S}}^m_{\rho,\delta,\rm{rg}}(K\times\R^p)$ is the set of all generalized symbols $a$ having a representative $(a_\eps)_\eps$ fulfilling the condition 
\beq
\label{cond_1}
\exists N\in\N\ \forall \alpha\in\N^p\, \forall\beta\in\N^n\qquad\qquad 
 |a_\eps|^{(m)}_{\rho,\delta,K,\alpha,\beta}=O(\eps^{-N})\ \text{as $\eps\to 0$},
\eeq
the space $\Neg_{\rho,\delta}(\Om\times\R^p)$ is the set $\Neg_{S^m_{\rho,\delta}(\Om\times\R^p)}$ of negligible nets, $\mathcal{S}_{\rho,\delta,\ssc}(\Om\times\R^p)$ is the set of all nets $(a_\eps)_\eps\in S^m_{\rho,\delta}(\Om\times\R^p)^{(0,1]}$ such that 
\beq
\label{cond_2}
\forall K\Subset\Om\, \exists (\omega_\eps)_\eps\in\Pi_\ssc\, \forall\alpha\in\N^p\, \forall\beta\in\N^n\qquad\qquad  |a_\eps|^{(m)}_{\rho,\delta,K,\alpha,\beta}=O(\omega_\eps)
\eeq
and finally $\wt{\mathcal{S}}^{-\infty}_{\rm{rg}}(K\times\R^p)$ is the space of all $a\in\wt{\mathcal{S}}^{-\infty}(\Om\times\R^p)$ for which there exists a representative $(a_\eps)_\eps$ with the property
\beq
\label{cond_3}
\exists N\in\N\, \forall m\in\R\, \forall \alpha\in\N^p\, \forall\beta\in\N^n\qquad\qquad  |a_\eps|^{(m)}_{K,\alpha,\beta}=O(\eps^{-N}).
\eeq
By construction $\wt{\mathcal{S}}^m_{\rho,\delta,\ssc}(\Om\times\R^p)\subseteq\wt{\mathcal{S}}^m_{\rho,\delta,\rm{rg}}(\Om\times\R^p)$ and the generalized symbols of $\wt{\mathcal{S}}^{-\infty}(\Om\times\R^p)$ and $\wt{\mathcal{S}}^{-\infty}_{\rm{rg}}(\Om\times\R^p)$ can be regarded as elements of $\wt{\mathcal{S}}^m_{\rho,\delta}(\Om\times\R^p)$ and $\wt{\mathcal{S}}^m_{\rho,\delta,\rm{rg}}(\Om\times\R^p)$ respectively. In all the previous notations the absence of the subindex $(\rho,\delta)$ means $\rho=1$ and $\delta=0$.

When $a\in\wt{\mathcal{S}}^m_{\rho,\delta}(\Om\times\Om\times\R^n)$ and $u\in\Gc(\Om)$ the generalized oscillatory integral 
\[
\int_{\Om\times\R^n}e^{i(x-y)\xi}a(x,y,\xi)u(y)\, dy\, \dslash\xi:=\biggl[\biggl(\int_{\Om\times\R^n}e^{i(x-y)\xi}a_\eps(x,y,\xi)u_\eps(y)\, dy\,\dslash\xi\biggr)_\eps\biggr]
\]
defines a generalized function in $\G(\Om)$ (see \cite{Garetto:04th, GGO:03} for details on the theory of generalized oscillatory integrals). The $\wt{\C}$-linear continuous map
\[
A:\Gc(\Om)\to\G(\Om):u\to\int_{\Om\times\R^n}e^{i(x-y)\xi}a(x,y,\xi)u(y)\, dy\, \dslash\xi
\]
is called \emph{generalized pseudodifferential operator} of amplitude $a\in\wt{\mathcal{S}}^m_{\rho,\delta}(\Om\times\Om\times\R^n)$. In the next proposition we recall the mapping properties concerning generalized pseudodifferential operators which will be involved in the extension of the action of $A$ from $\Gc(\Om)$ to $\LL(\G(\Om),\wt{\C})$. For the corresponding proofs the reader should refer to \cite[Chapter 4]{Garetto:04th} and \cite[Section 4]{GGO:03}.
\begin{proposition}
\label{prop_pseudo_1}
\leavevmode
\begin{itemize}
\item[(i)] Let $a\in\wt{\mathcal{S}}^m_{\rho,\delta}(\Om\times\Om\times\R^n)$. The corresponding pseudodifferential operator $A$ is a continuous map from $\Gc(\Om)$ to $\G(\Om)$.
\item[(ii)] If $a\in\wt{\mathcal{S}}^m_{\rho,\delta,\rm{rg}}(\Om\times\Om\times\R^n)$ then $A$ maps $\Gcinf(\Om)$ continuously into $\Ginf(\Om)$.
\item[(iii)] If $a\in\wt{\mathcal{S}}^{-\infty}_{\rm{rg}}(\Om\times\Om\times\R^n)$ then $A$ maps $\Gc(\Om)$ continuously into $\Ginf(\Om)$.
\end{itemize}
\end{proposition}
Clearly Proposition \ref{prop_pseudo_1} can be stated for the formal transposed of $A$. ${\,}^tA$ is the pseudodifferential operator defined by $${\ }^tAu(x)=\int_{\Om\times\R^n}e^{i(x-y)\xi}a(y,x,-\xi)u(y)\,dy\,\dslash\xi .$$
The functional $k_A\in\LL(\Gc(\Om\times\Om),\wt{\C})$ given by
\[
k_A(u):=\int_{\Om}A(u(x,\cdot))\, dx
\]
is called \emph{kernel} of $A$. As shown in \cite{Garetto:04th} it can be written as the oscillatory integral $$\int_{\Om\times\Om\times\R^n}e^{i(x-y)\xi}a(x,y,\xi)u(x,y)\, dx\, dy\, \dslash\xi$$ and fulfills the property 
\[
k_A(v\otimes u)=\int_{\Om}Au(x)v(x)\, dx =\int_\Om u(y){\,}^tAu(y)\, dy
\]
for all $u,v\in\Gc(\Om)$ ($v\otimes u:=(v_\eps(x)u_\eps(y))_\eps+\Neg(\Om\times\Om)$). We say that the pseudodifferential operator $A$ is \emph{properly supported} if the support of its kernel is a proper subset of $\Om\times\Om$. Proposition 4.3.18 in \cite{Garetto:04th} proves that the following mapping properties hold for properly supported pseudodifferential operators.
\begin{proposition}
\label{prop_ext_pseudo}
If $A$ is a properly supported pseudodifferential operator with amplitude $a\in\wt{\mathcal{S}}^m_{\rho,\delta}(\Om\times\Om\times\R^n)$ then
\begin{itemize}
\item[(i)] $A$ maps $\Gc(\Om)$ continuously into itself,
\item[(ii)] $A$ can be uniquely extended to a $\wt{\C}$-linear continuous map from $\G(\Om)$ into $\G(\Om)$ such that for all $u\in\G(\Om)$ and $v\in\Gc(\Om)$,
\begin{equation}
\label{inttrans}
\int_\Omega Au(x)v(x)dx=\int_\Omega u(y){\ }^tAv(y)dy.
\end{equation}
\end{itemize}
\ \ If $a\in\wt{\mathcal{S}}^m_{\rho,\delta,{\rm{rg}}}(\Om\times\Om\times\R^n)$ then
\begin{itemize}
\item[(iii)] $A$ maps $\Gcinf(\Om)$ continuously into itself,
\item[(iv)] the extension defined above maps $\Ginf(\Om)$ continuously into itself.
\end{itemize}
The same results hold with ${\ }^tA$ in place of $A$.
\end{proposition}
Concluding, let us consider the expression 
\beq
\label{integral_expr}
Ru(x)=\int_{\Om}k(x,y)u(y)\, dy,
\eeq
where $k\in\G(\Om\times\Om)$ and $u\in\Gc(\Om)$. By Proposition 1.2.25 in \cite{Garetto:04th} \eqref{integral_expr} defines a continuous $\wt{\C}$-linear operator $R:\Gc(\Om)\to\G(\Om)$. Note that $k$ is uniquely determined by \eqref{integral_expr} as an element of $\G(\Om\times\Om)$. For this reason, we may call it {\em the} kernel of $R$, adopt the notation $k_R$, and we may call $R$ an \emph{operator with generalized kernel}. When $k_R\in\Ginf(\Om\times\Om)$ then $R$ maps $\Gc(\Om)$ continuously into $\Ginf(\Om)$ and we use the expression \emph{operator with regular generalized kernel}. A simple adaptation of the reasoning of \cite[Proposition 4.3.13]{Garetto:04th} and \cite[Proposition 4.12]{GGO:03} yields the following characterizations.
\begin{proposition}
\label{prop_charac_infty}
\leavevmode
\begin{trivlist}
\item[(i)] $R$ is an operator with generalized kernel if and only if it is a pseudodifferential operator with amplitude in $\wt{\mathcal{S}}^{-\infty}(\Om\times\Om\times\R^n)$.
\item[(ii)] $R$ is an operator with regular generalized kernel if and only if it is a pseudodifferential operator with amplitude in $\wt{\mathcal{S}}^{-\infty}_{\rm{rg}}(\Om\times\Om\times\R^n)$.
\end{trivlist}
\end{proposition}
\begin{proof}
We leave to the reader to check that in both $(i)$ and $(ii)$ the amplitude $r$ is given by $$e^{i(x-y)\xi}k_R(x,y)\chi(\xi):=(e^{i(x-y)\xi}k_{R,\eps}(x,y)\chi(\xi))_\eps+\Neg^{-\infty}(\Om\times\Om\times\R^n),$$ where $\chi$ is a cut-off function in $\Cinfc(\R^n)$ with $\int\chi(\xi)\, \dslash\xi=1$.
\end{proof}
\begin{remark}
\label{rem_basic_pseudo}
According to the language of Section \ref{section_basic} all the operators considered so far are basic. More precisely, the generalized pseudodifferential operators are basic elements of the space $\LL(\Gc(\Om),\G(\Om))$ and their kernels are basic functionals in $\LL(\Gc(\Om\times\Om),\wt{\C})$.
\end{remark}
\begin{definition}
\label{def_pseudo_extension}
Let $A$ be a pseudodifferential operator with amplitude $a\in\wt{\mathcal{S}}^m_{\rho,\delta}(\Om\times\Om\times\R^n)$. We extend the action of $A$ to the dual $\LL(\G(\Om),\wt{\C})$ as
\beq
\label{ext_formula}
AT(u):=T({\, }^tAu),\qquad\qquad\qquad u\in\Gc(\Om).
\eeq
\end{definition}
The fact that \eqref{ext_formula} extends the original definition of $A$ on $\Gc(\Om)$ is due to the equality \eqref{inttrans}.
\begin{proposition}
\label{prop_pseudo_extension}
\leavevmode
\begin{itemize}
\item[(i)] The operator $A$ defined in \eqref{ext_formula} maps $\LL(\G(\Om),\wt{\C})$ continuously into the dual $\LL(\Gc(\Om),\wt{\C})$.
\item[(ii)] If $A$ is properly supported then it maps $\LL(\G(\Om),\wt{\C})$ and $\LL(\Gc(\Om),\wt{\C})$ into themselves respectively and with continuity. 
\item[(iii)] If $a\in\wt{\mathcal{S}}^m_{\rho,\delta,{\rm{rg}}}(\Om\times\Om\times\R^n)$ then $A$ maps $\LL(\Ginf(\Om),\wt{\C})$ continuously into $\LL(\Gcinf(\Om),\wt{\C})$ and when it is properly supported the duals $\LL(\Ginf(\Om),\wt{\C})$ and $\LL(\Gcinf(\Om),\wt{\C})$ are mapped into themselves respectively with continuity.
\end{itemize}
\end{proposition} 
\begin{proof}
$(i)$ By composition of continuous maps $AT$ is an element of the dual $\LL(\Gc(\Om),\wt{\C})$ when $T\in\LL(\G(\Om),\wt{\C})$. Since the map ${\,}^tA:\Gc(\Om)\to\G(\Om)$ is continuous the image ${\,}^tA(B)$ of a bounded subset $B$ of $\Gc(\Om)$ is bounded in $\G(\Om)$. Hence, from 
\beq
\label{cont_A_1}
\sup_{u\in B}|AT(u)|_\esp = \sup_{v\in {\,}^tA(B)}|Tv|_\esp
\eeq
we have that $A$ is continuous from $\LL(\G(\Om),\wt{\C})$ to $\LL(\Gc(\Om),\wt{\C})$.

$(ii)$ By the assertions $(i)$ and $(ii)$ of Proposition \ref{prop_ext_pseudo} if $A$ is properly supported then ${\,}^tA$ maps $\Gc(\Om)$ continuously into itself and can be extended to a $\wt{\C}$-linear continuous map on $\G(\Om)$. It follows that $A$ maps the duals $\LL(\G(\Om),\wt{\C})$ and $\LL(\Gc(\Om),\wt{\C})$ into themselves respectively and with continuity.

$(iii)$ Assume that $a\in\wt{\mathcal{S}}^m_{\rho,\delta,{\rm{rg}}}(\Om\times\Om\times\R^n)$. By Proposition \ref{prop_pseudo_1}$(ii)$ the operator ${\,}^tA$ is continuous from $\Gcinf(\Om)$ to $\Ginf(\Om)$. This means that if $T\in\LL(\Ginf(\Om),\wt{\C})$ then $AT\in\LL(\Gcinf(\Om),\wt{\C})$. Since ${\,}^tA(B)$ is bounded in $\Ginf(\Om)$ when $B$ is a bounded subset of $\Gcinf(\Om)$ we have that the equality \eqref{cont_A_1} holds and proves the continuity of $A$ from $\LL(\Ginf(\Om),\wt{\C})$ to $\LL(\Gcinf(\Om),\wt{\C})$. By combining $(iii)$ and $(iv)$ in Proposition \ref{prop_ext_pseudo} with the definition given in \eqref{ext_formula} the proof of the assertion is complete.
\end{proof} 
The action of a generalized pseudodifferential operator on a basic functional gives a functional which is still basic in all the statements above. 

Due to Proposition \ref{prop_charac_infty} the action of an operator with generalized kernel on a functional in $\LL(\G(\Om),\wt{\C})$ can be seen as the action of a pseudodifferential operator with generalized symbol of order $-\infty$ on a functional. Interesting mapping properties are obtained on the spaces $\Lb(\G(\Om),\wt{\C})$ and $\Lb(\Gc(\Om),\wt{\C})$ of basic functionals. 
\begin{proposition}
\label{prop_reg_func}
Let $R$ be an operator with generalized kernel $k_R\in\G(\Om\times\Om)$.
\begin{itemize}
\item[(i)] $R$ maps $\Lb(\G(\Om),\wt{\C})$ into $\G(\Om)$ with continuity.
\item[(ii)] If $R$ is properly supported then it is a continuous map from $\Lb(\G(\Om),\wt{\C})$ to $\Gc(\Om)$ and from $\Lb(\Gc(\Om),\wt{\C})$ to $\G(\Om)$.
\item[(iii)] If $k_R\in\Ginf(\Om\times\Om)$ then $(i)$ and $(ii)$ hold for the same spaces of functionals and with $\Gcinf(\Om)$ and $\Ginf(\Om)$ in place of $\Gc(\Om)$ and $\G(\Om)$ respectively.
\end{itemize}
\end{proposition}
\begin{proof}
$(i)$ By Definition \ref{def_pseudo_extension} we know that for $T\in\LL(\G(\Om),\wt{\C})$ and $u\in\Gc(\Om)$  
\[
RT(u)=T({\,}^tRu)=T\biggl(\int_\Om k_R(x,\cdot)u(x)\, dx\biggr),
\]
where $\int_\Om k_R(x,\cdot)u(x)\, dx\in\G(\Om)$. Since $T$ is basic, from Proposition \ref{prop_func_int}$(i)$ we have that 
\beq
\label{equality_TR}
T\biggl(\int_\Om k_R(x,\cdot)u(x)\, dx\biggr) = \int_\Om T(k_R(x,\cdot))u(x)\, dx
\eeq
and then 
\[
RT(u)=\int_\Om T(k_R(x,\cdot))u(x)\, dx.
\]
As shown in Proposition \ref{prop_action_func}$(v)$, $T(k_R(x,\cdot))$ is a generalized function in $\G(\Om)$. This proves that $RT\in\G(\Om)$. Note that for $K\Subset\Om$ the subset $B_{K,j,R}:=\{v:=\partial^\alpha_x k_R(x,\cdot)\}_{|\alpha|\le j, x\in K}$ of $\G(\Om)$ is bounded. Easy computations at the level of representatives lead to the inequality
\beq
\label{inequality_2}
\mP_{K,j}(RT)\le\sup_{v\in B_{K,j,R}}|Tv|_\esp
\eeq
which proves the continuity of $R$ from $\Lb(\G(\Om),\wt{\C})$ to $\G(\Om)$.

$(ii)$ If $R$ is properly supported then its kernel $k_R$ has proper support and by Proposition \ref{prop_action_func}$(vii)$ the generalized function $T(k_R(x,\cdot))\in\Gc(\Om)$. An arguing analogous to the one employed in \eqref{inequality_2} yields that $R$ is continuous from $\Lb(\G(\Om),\wt{\C})$ to $\Gc(\Om)$. Assume now that $T$ is a basic functional of $\LL(\Gc(\Om),\wt{\C})$. Since $k_R$ has proper support from Proposition \ref{prop_action_func}$(iii)$ we have that $T(k_R(x,\cdot))\in\G(\Om)$ and by Proposition \ref{prop_func_int}$(iv)$ we conclude that \eqref{equality_TR} is valid for all $u\in\Gc(\Om)$. Thus, $R$ maps $\Lb(\Gc(\Om),\wt{\C})$ into $\G(\Om)$. Since the set $C_{K,j,R}:=\{v:=\partial^\alpha_x k_R(x,\cdot)\chi(x,\cdot)\}_{|\alpha|\le j, x\in K}$, where $\chi$ is a proper function identically $1$ in a neighborhood of $\supp\, k_R$ is contained in $\Gc(\Om)$ and bounded there, the inequality
\[
\mP_{K,j}(RT)\le \sup_{v\in C_{K,j,R}}|Tv|_\esp
\]
implies the continuity of $R:\Lb(\Gc(\Om),\wt{\C})\to\G(\Om)$.

$(iii)$ Finally we suppose that $k_R\in\Ginf(\Om\times\Om)$. The desired mapping properties are a consequence of the assertions $(vi)$, $(viii)$ and $(iv)$ in Proposition \ref{prop_action_func}. For what concerns the continuity, an investigation at the representatives'level shows that if $T\in\Lb(\G(\Om),\wt{\C})$ then
\[
\mP_{\Ginf(K)}(RT)\le\sup_{v\in B_{K,R}}|Tv|_\esp,
\]
where, for $\psi$ cut-off function in $\Cinfc(\Om)$ identically 1 in a neighborhood of $\supp\, T$, the set $B_{K,R}:=\{v:=\partial^\alpha_x k_R(x,\cdot)\psi(\cdot)\}_{\alpha\in\N^n,x\in K}\subseteq\G(\Om)$ is bounded. Analogously if $\supp\, k_R$ is proper then $R$ is continuous from $\Lb(\G(\Om),\wt{\C})$ to $\Gcinf(\Om)$. Since the set $C_{K,R}:=\{v:=\partial^\alpha_x k_R(x,\cdot)\chi(x,\cdot)\}_{\alpha\in\N^n,x\in K}$ is bounded in $\Gc(\Om)$, the inequality
\[
\mP_{\Ginf(K)}(RT)\le\sup_{v\in C_{K,R}}|Tv|_\esp,
\]
valid for $T$ in $\Lb(\Gc(\Om),\wt{\C})$ entails the continuity of $R$ from $\Lb(\Gc(\Om),\wt{\C})$ to $\Ginf(\Om)$.
\end{proof}
Proposition 4.10 in \cite{GGO:03} has a natural version in the dual context.
\begin{proposition}
\label{prop_W_W'}
Let $A$ be a pseudodifferential operator with amplitude $a\in\wt{\mathcal{S}}^m_{\rho,\delta}(\Om\times\Om\times\R^n)$.
\begin{itemize}
\item[(i)] $k_A\in\G(\Om\times\Om\setminus\Delta)$ where $\Delta$ is the diagonal of $\Om\times\Om$ and for $W$ and $W'$ open subsets of $\Om$ with $W\times W'\subseteq\Om\times\Om\setminus\Delta$, the equality
\[
Au|_W(x)=\int_{W'}k_A(x,y)u(y)\, dy
\]
holds for all $u\in\Gc(W')$. Moreover, for all $u\in\Gc(W)$,
\beq
\label{transposed_eq}
{\ }^tAu|_{W'}(x)=\int_W k_A(y,x)u(y)\, dy.
\eeq
\item[(ii)] If $T$ is a basic functional of $\LL(\G(\Om),\wt{\C})$ with $\supp\, T\subseteq W'$ then $AT|_W\in\G(W)$ and 
\[
AT|_W(u)=\int_W T(k_A(x,\cdot))u(x)\, dx
\]
for all $u\in\Gc(W)$.
\item[(iii)] If $a\in\wt{\mathcal{S}}^m_{\rho,\delta,{\rm{rg}}}(\Om\times\Om\times\R^n)$ then $k_A\in\Ginf(\Om\times\Om\setminus\Delta)$ and for all $T$ basic functional of $\LL(\G(\Om),\wt{\C})$ with $\supp\, T\subseteq W'$ the restriction $AT|_W$ belongs to $\Ginf(W)$.
\end{itemize}
\end{proposition}
\begin{proof}
The first assertion is proven in \cite[Proposition 4.10]{GGO:03}. In particular the equality \eqref{transposed_eq} is due to \eqref{inttrans}. Let now $T$ be a basic functional of $\LL(\G(\Om),\wt{\C})$ with $\supp\, T\subseteq W'$ and $\psi$ be a cut-off function in $\Cinfc(W')$ identically $1$ in a neighborhood of $\supp\, T$. By \eqref{transposed_eq} for all $u\in\Gc(W)$ we can write
\[
AT|_W(u)=T({\,}^tA u)=T(({\,}^tAu)\psi)=T\biggl(\int_W k_A(y,\cdot)u(y)\, dy\,\psi(\cdot)\biggr).
\] 
From Proposition \ref{prop_func_int}$(i)$ and Proposition \ref{prop_action_func}$(v)$ the functional $T$ goes under the integral sign and we obtain that
\beq
\label{eq_AT_W}
AT|_W(u)=\int_W T(k_A(y,\cdot)\psi(\cdot))u(y)\, dy = \int_W T(k_A(x,\cdot))u(x)\, dx,
\eeq
where $T(k_A(x,\cdot))\in\G(W)$. Finally, when $a\in\wt{\mathcal{S}}^m_{\rho,\delta,{\rm{rg}}}(\Om\times\Om\times\R^n)$ then $k_A\in\Ginf(\Om\times\Om\setminus\Delta)$ (see \cite[Proposition 4.10$(iv)$]{GGO:03}) and combining Proposition \ref{prop_func_int}$(i)$ with Proposition \ref{prop_action_func}$(vi)$ we have that $T(k_A(x,\cdot))\in\Ginf(W)$. In view of \eqref{eq_AT_W} this means that $AT|_W\in\Ginf(W)$.
\end{proof}
Before proceeding we observe that since the dual $\LL(\Gc(\Om),\wt{\C})$ contains both the Colombeau algebras $\Ginf(\Om)$ and $\G(\Om)$ by continuous embedding and $\Om\to\LL(\Gc(\Om),\wt{\C})$ is a sheaf, it is meaningful to look for the regions where a functional is a Colombeau function or a $\Ginf$-Colombeau function. The regularity with respect to $\G(\Om)$ or with respect to $\Ginf(\Om)$ is measured by the following notions of \emph{$\G$-singular support} and \emph{$\Ginf$-singular support}. 
\begin{definition}
\label{def_sing_supports}
The $\G$-singular support of $T\in\LL(\Gc(\Om),\wt{\C})$ ($\singsupp_\G\, T$) is the complement of the set of all points $x\in\Om$ such that the restriction of $T$ to some neighborhood $V$ of $x$ belongs to $\G(V)$.\\
The $\Ginf$-singular support of $T\in\LL(\Gc(\Om),\wt{\C})$ ($\singsupp_{\Ginf}\, T$) is the complement of the set of all points $x\in\Om$ such that the restriction of $T$ to some neighborhood $V$ of $x$ belongs to $\Ginf(V)$. 
\end{definition}
By definition it is clear that $\singsupp_\G\, T$ and $\singsupp_{\Ginf}\, T$ are both closed subsets of $\Om$ and that the inclusion $\singsupp_\G\, T\subseteq\singsupp_{\Ginf}\, T$ holds. Clearly the $\Ginf$-singular support extends the usual notion of generalized singular support of a Colombeau function in $\G(\Om)$ (\cite[Section 2]{GGO:03}) to the dual $\LL(\Gc(\Om),\wt{\C})$.

The pseudolocality-property proved in \cite{GGO:03} to be valid for a pseudodifferential operator $A:\Gc(\Om)\to\G(\Om)$ with regular amplitude can now be stated for the extension to the dual $\LL(\G(\Om),\wt{\C})$. The two different ways of measuring the regularity of a functional in $\LL(\G(\Om),\wt{\C})$ considered above give a new and more elaborated aspect to the following result of pseudolocality.
\begin{theorem}
\label{theorem_pseudolocality_dual} 
Let $T$ be a basic functional in $\LL(\G(\Om),\wt{\C})$.
\begin{itemize}
\item[(i)] If $A$ is a pseudodifferential operator with amplitude in $\wt{\mathcal{S}}^m_{\rho,\delta}(\Om\times\Om\times\R^n)$ then
\beq
\label{pseudolocality_1}
\singsupp_\G\,AT\subseteq\singsupp_\G\, T.
\eeq
\item[(ii)] If $A$ is a pseudodifferential operator with amplitude in $\wt{\mathcal{S}}^m_{\rho,\delta,{\rm{rg}}}(\Om\times\Om\times\R^n)$ then
\beq
\label{pseudolocality_2}
\singsupp_{\Ginf}\,AT\subseteq\singsupp_{\Ginf}\, T.
\eeq
\end{itemize}
\end{theorem}
It is clear that \eqref{pseudolocality_1} and \eqref{pseudolocality_2} can be written for basic functionals in the dual $\LL(\Gc(\Om),\wt{\C})$ when $A$ is properly supported.
\begin{proof}
$(i)$ Given a basic functional $T\in\LL(\G(\Om),\wt{\C})$ we consider an open neighborhood $V$ of $\singsupp_\G\, T$ and $\psi\in\Cinfc(V)$ identically $1$ in a neighborhood of $\singsupp_\G\, T$. We write $T=\psi T+(1-\psi)T$. $\psi T$ is a basic functional in $\LL(\G(\Om),\wt{\C})$ and by definition of $\G$-singular support we know that $(1-\psi)T\in\Gc(\Om)$. From Proposition \ref{prop_pseudo_1}$(i)$ we have that $A((1-\psi)T)\in\G(\Om)$, hence our assertion becomes
\beq
\label{asser_A}
\singsupp_\G\, A(\psi T)\subseteq\singsupp_\G\, T.
\eeq
We will prove \eqref{asser_A} by preliminary showing that 
\beq
\label{asser_B}
\singsupp_\G\, AT\subseteq \supp\, T.
\eeq
Let $x_0\in\Om\setminus\supp\, T$ and $W$ and $W'$ be open neighborhoods of $x_0$ and $\supp\,T$ respectively such that $W\times W'\subseteq\Om\times\Om\setminus\Delta$. By Proposition \ref{prop_W_W'} it follows that $AT|_W\in\G(W)$ and therefore \eqref{asser_B} is proven. Replacing $T$ with $\psi T$ in \eqref{asser_B} we get $\singsupp_\G\, A(\psi T)\subseteq\supp\, \psi T\subseteq V$ and since $V$ is arbitrary the proof is complete.

$(ii)$ Assume now that $V$ is an open neighborhood of $\singsupp_{\Ginf}\, T$. By definition of $\Ginf$-singular support, $(1-\psi)T\in\Gcinf(\Om)$ and from Proposition \ref{prop_pseudo_1}$(ii)$, $A((1-\psi)T)\in\Ginf(\Om)$. Thus, our assertion becomes $\singsupp_{\Ginf} AT\subseteq \supp\, T$. This can be proven as above making use of Proposition \ref{prop_W_W'}$(iii)$.
\end{proof}
A parametrix construction, based on the symbolic calculus for generalized pseudodifferential operators developed in \cite{GGO:03}, can be provided for pseudodifferential operators whose generalized symbols satisfy suitable assumptions of hypoellipticity. In the sequel, we slightly simplify the notion of generalized hypoelliptic symbol introduced in \cite{GGO:03} and within the dual context we state the theorem on the existence of a parametrix and the following result of Colombeau regularity. For technical reasons (see \cite[Section 5]{GGO:03}) we consider generalized symbols $a$ whose representing nets $(a_\eps)_\eps$ fulfill the characterizing seminorms estimates for all values of $\eps$ in the interval $(0,1]$. In this way we define the subspaces $\underline{\wt{\mathcal{S}}}^{\,m}_{\,\rho,\delta}(\Om\times\R^p)$, $\underline{\wt{\mathcal{S}}}^{\,m}_{\,\rho,\delta,{\rm{rg}}}(\Om\times\R^p)$, $\underline{\wt{\mathcal{S}}}^{\,m}_{\,\rho,\delta,\ssc}(\Om\times\R^p)$, $\underline{\wt{\mathcal{S}}}^{-\infty}(\Om\times\R^p)$ and $\underline{\wt{\mathcal{S}}}^{-\infty}_{\,\rm{rg}}(\Om\times\R^p)$ of ${\wt{\mathcal{S}}}^{\,m}_{\rho,\delta}(\Om\times\R^p)$, ${\wt{\mathcal{S}}}^{\,m}_{\rho,\delta,{\rm{rg}}}(\Om\times\R^p)$, ${\wt{\mathcal{S}}}^{\,m}_{\rho,\delta,\ssc}(\Om\times\R^p)$, ${\wt{\mathcal{S}}}^{-\infty}(\Om\times\R^p)$ and ${\wt{\mathcal{S}}}^{-\infty}_{\rm{rg}}(\Om\times\R^p)$ respectively.
\begin{definition}
\label{def_gen_hyp}
Let $m,l,\rho,\delta$ be real numbers with $l\le m$ and $0\le\delta<\rho\le 1$. We say that $a\in\underline{\wt{\mathcal{S}}}^{\,m}_{\,\rho,\delta,{\rm{rg}}}(\Om\times\R^n)$ is a generalized hypoelliptic symbol  of order $(m,l)$ and type $(\rho,\delta)$ if it has a representative $(a_\eps)_\eps$ fulfilling the following condition: 
for all $K\Subset\Omega$ there exists a strongly positive slow scale net
$(r_{K,\eps})_\eps$, a net $(\omega_{1,K,\eps})_\eps$, $\omega_{1,K,\eps}\ge C_K\eps^{s_K}$ on the
interval $(0,1]$ for certain constants $C_K>0$, $s_K\in\mathbb{R}$, and slow scale nets
$(\omega_{2,K,\alpha,\beta,\eps})_\eps$, such that for all $x\in K$, for $|\xi|\ge r_{K,\eps}$,
for all $\eps\in(0,1]$,
\begin{equation}
\label{hyp1}
|a_\eps(x,\xi)|\ge \omega_{1,K,\eps}\langle\xi\rangle^l
\end{equation}
and
\begin{equation}
\label{hyp2}
|\partial^\alpha_\xi\partial^\beta_x a_\eps(x,\xi)|\le \omega_{2,K,\alpha,\beta,\eps}|a_\eps(x,\xi)|\langle\xi\rangle^{-\rho|\alpha|+\delta|\beta|}.
\end{equation}
for all $(\alpha,\beta)\neq(0,0)$.
\end{definition} 
A generalized symbol $a$ satisfying Definition \ref{def_gen_hyp} with $l=m$ is said to be \emph{elliptic}. 

Theorem 6.8 in \cite{GGO:03} proves that when $A$ is a pseudodifferential operator with generalized hypoelliptic symbol of order $(m,l)$ and type $(\rho,\delta)$ then there exists a properly supported pseudodifferential operator $P$ with symbol in $\underline{\wt{\mathcal{S}}}^{-l}_{\rho,\delta,{\rm{rg}}}(\Om\times\R^n)$ such that for all $u\in\Gc(\Om)$,
\beq
\label{parametrix_eq}
\begin{split}
PAu&=u+Ru,\\
APu&=u+Su,
\end{split}
\eeq
where $R$ and $S$ are operators with regular generalized kernel. Note that if $A$ properly supported then $R$ and $S$ are properly supported operators themselves and the equalities in \eqref{parametrix_eq} hold for all $u\in\G(\Om)$. Moreover, by definition of the extension of a pseudodifferential operator to the dual $\LL(\Gc(\Om),\wt{\C})$ we can replace $u$ with $T\in\LL(\Gc(\Om),\wt{\C})$ in \eqref{parametrix_eq}. In this case the equalities have to be read in $\LL(\Gc(\Om),\wt{\C})$. 
\begin{theorem}
\label{theo_reg_dual}
Let $A$ be a properly supported pseudodifferential operator with generalized hypoelliptic symbol. Then, for every $T$ basic functional in the dual $\LL(\Gc(\Om),\wt{\C})$,
\beq
\label{pseudoloc_A}
\singsupp_\G\, AT= \singsupp_\G\, T
\eeq
and
\beq
\label{pseudoloc_B}
\singsupp_{\Ginf}\, AT= \singsupp_{\Ginf}\, T.
\eeq
\end{theorem}
\begin{proof}
The inclusions $\singsupp_\G\, AT\subseteq \singsupp_\G\, T$ and $\singsupp_{\Ginf}\, AT\subseteq \singsupp_{\Ginf}\, T$ are clear from the pseudolocality-property of $A$ (Theorem \ref{theorem_pseudolocality_dual}). Taking a parametrix $P$ we can write $T$ as $PAT-RT$ and by the assertions $(ii)$ and $(iii)$ of Proposition \ref{prop_reg_func}, $RT$ belongs to $\Ginf(\Om)$. Hence, $\singsupp_\G\, T=\singsupp_\G\, PAT$ and $\singsupp_{\Ginf}\, T=\singsupp_{\Ginf}\, PAT$. At this point since $AT$ is a basic functional of $\LL(\Gc(\Om),\wt{\C})$ the pseudolocality-property of $P$ allows to conclude the proof.
\end{proof}

\section{$\G$-wave front set and $\Ginf$-wave front set of a functional in $\LL(\Gc(\Om),\wt{\C})$}
\label{section_WF}
Microlocal analysis in Colombeau algebras of generalized functions as it has been initiated (in published form) in \cite{DPS:98, NPS:98} is a compatible extension of its distribution theoretic analogue to an unrestricted differential-algebraic context. The classical H\"ormander definition of wave front set of a distribution $u$ (see \cite{Hoermander:71}) makes use of the notion of micro-ellipticity and consists in the intersection of the characteristic sets (i.e. region of non-ellipticity) of those pseudodifferential operators which map $u$ in a $\Cinf$-function. A characterization of $\WF u$ is given in terms of direct estimates of the Fourier transform of $u$, after multiplication by a suitable cut-off function. 

The \emph{generalized wave front set of $u\in\G(\Om)$} (or \emph{$\Ginf$-wave front set of $u$} denoted by $\WF_{\Ginf}(u)$) is defined by translating the Fourier transform-characterization of the distributional wave front set into the language of representatives of generalized functions and replacing the $\Cinf$-regularity with the $\Ginf$-regularity. This sort of ``elementary'' approach to the wave front set is a natural definition in the Colombeau framework. 

The theory of generalized pseudodifferential operators established in \cite{GGO:03} and extended to the dual space $\LL(\Gc(\Om),\wt{\C})$ in the previous section, has suggested a ``pseudodifferential-characterization'' of the $\Ginf$-wave front set of a Colombeau generalized function. This has been provided in \cite{GH:05} by making use of pseudodifferential operators with slow scale symbols and introducing a sufficiently strong notion of micro-ellipticity. In the sequel we give an essential overview of the concepts in \cite{GH:05} which will be employed in this section more frequently, referring to \cite{GH:05} for the proofs of the main results and for further explanations.
\begin{definition}
\label{def_micro_ellipticity}
Let $a \in \wt{{{\mathcal{S}}}}_{\rho,\delta,\ssc}^m(\Om\times\R^n)$ and $(x_0,\xi_0) \in \CO{\Om}$. We say that $a$ is slow scale micro-elliptic at $(x_0,\xi_0)$ if it has a representative $(a_\eps)_\eps$ satisfying the following: there is a relatively compact open neighborhood $U$ of $x_0$, a conic neighborhood $\Gamma$ of $\xi_0$, and $(r_\eps)_\eps , (s_\eps)_\eps$ in $\Pi_\ssc$ such that
\beq
\label{estimate_below}
 | a_\eps(x,\xi)| \ge \frac{1}{s_\eps} \lara{\xi}^m\qquad\qquad (x,\xi)\in U\times\Gamma,\, |\xi| \ge r_\eps,\, \eps \in (0,1].
\eeq
We denote by $\Ellsc(a)$ the set of all $(x_0,\xi_0) \in \CO{\Om}$ where $a$ is slow scale micro-elliptic.\\
If there exists $(a_\eps)_\eps \in\ a$ such that \eqref{estimate_below} holds at all points in $\CO{\Om}$ then the symbol $a$ is called \emph{slow scale elliptic}.
\end{definition}
\begin{remark}
\label{rem_slow_ell}
\leavevmode
\begin{trivlist}
\item[(i)] Note that in the definition of the set $\Ellsc(a)$ makes no difference to require that the estimate from below in \eqref{estimate_below} holds for all $\eps\in(0,1]$ or in a smaller interval $(0,\eta]$. Indeed, assume that \eqref{estimate_below} holds for some representative $(a_\eps)_\eps$ of $a$ when $\eps$ is smaller of a certain $\eta\in(0,1]$ and take $b\in S^m_{\rho,\delta}(\Om\times\R^n)$ such that $|b(x,\xi)|\ge\lara{\xi}^m$ on $U\times\Gamma$. It is not restrictive to suppose that the representative $(a_\eps)_\eps$ is identically $0$ when $(x,\xi)\in U\times\Gamma$ and $\eps\in(\eta,1]$. Let $(\omega_\eps)_\eps$ be a net in $\R^{(0,1]}$ defined as follows: $\omega_\eps=1/s_\eps$ for $\eps\in(\eta,1]$, $\omega_\eps=0$ for $\eps\in(0,\eta]$. It is clear that $a'_\eps:=a_\eps+\omega_\eps b$ is another representative of $a$. Moreover, by construction, $|a'_\eps(x,\xi)|\ge s_\eps^{-1}\lara{\xi}^m$ when $(x,\xi)\in U\times\Ga$, $|\xi|\ge r_\eps$, $\eps\in(0,\eta]$ and $|a'_\eps(x,\xi)|=\omega_\eps|b(x,\xi)|\ge s_\eps^{-1}\lara{\xi}^m$ when $(x,\xi)\in U\times\Ga$, $|\xi|\ge r_\eps$, $\eps\in(\eta,1]$.
\item[(ii)] Any symbol $a \in \wt{\underline{{\mathcal{S}}}}_{\rho,\delta,\ssc}^m(\Om\times\R^n)$ which is slow scale micro-elliptic at $(x_0,\xi_0)$ fulfills the stronger hypoellipticity estimates of Definition \ref{def_gen_hyp} and it is stable under lower order (slow scale) perturbations \cite[Proposition 2.3]{GH:05}. More precisely if $(a_\eps)_\eps\in\mathcal{S}^m_{\rho,\delta,\ssc}(\Om\times\R^n)$ satisfy \eqref{estimate_below} in $U\times\Ga\ni(x_0,\xi_0)$ then
\begin{trivlist}
\item[-] for all $\alpha,\beta\in\N^n$ there exist $(\lambda_\eps)_\eps\in\Pi_\ssc$ and $\eta\in(0,1]$ such that
\[
|\partial^\alpha_\xi\partial^\beta_x a_\eps(x,\xi)|\le \lambda_\eps |a_\eps(x,\xi)|\lara{\xi}^{-\rho|\alpha|+\delta|\beta|},\qquad (x,\xi)\in U\times\Ga,\, |\xi|\ge r_\eps,\, \eps\in(0,\eta];
\]
\item[-] for all $(b_\eps)_\eps \in \mathcal{S}_{\rho,\delta,\ssc}^{m'}(\Om\times\R^n)$, $m'<m$, there exist $(r'_\eps)_\eps$, $(s'_\eps)_\eps \in \Pi_\ssc$ and $\eta\in(0,1]$ such that 
\[
|a_\eps(x,\xi)+b_\eps(x,\xi)| \ge \frac{1}{s'_\eps}\lara{\xi}^m\qquad (x,\xi) \in U\times\Gamma,\, |\xi|\ge r'_\eps,\, \eps\in(0,\eta].
\]
\end{trivlist}
\end{trivlist}
\end{remark}
As in \cite{GH:05} we choose $\Oprop{m}(\Om)$ for denoting sets of all properly supported operators $a(x,D)$ with symbol in $\Syscu^{m}(\Om\times\R^n)$ and given $u\in\G(\Om)$ we define the set
\beq
\label{set_Wsc_G}
{\rm{W}}^{\ssc}_{\Ginf}(u):=\hskip-5pt \bigcap_{\substack{a(x,D)\in\,\Oprop{0}(\Om)\\[0.1cm] a(x,D)u\, \in\, \Ginf(\Om)}}\hskip-5pt \compl{\Ellsc(a)}.
\eeq
Theorem 3.10 in \cite{GH:05} proves that for all $u\in\G(\Om)$,
\[
{\rm{W}}^{\ssc}_{\Ginf}(u)=\WF_{\Ginf}(u).
\]
Inspired by these results and aware of the fact that two kinds of regularity, with respect to $\G(\Om)$ and with respect to $\Ginf(\Om)$, coexist in the dual $\LL(\Gc(\Om),\wt{\C})$, in this section we define the \emph{$\G$-wave front set and the $\Ginf$-wave front set of $T\in\LL(\Gc(\Om),\wt{\C})$} and we provide a Fourier transform-characterization in case of basic functionals.

\subsection{Definition and basic properties of the generalized wave front sets $\WF_\G(T)$ and $\WF_{\Ginf}(T)$}
\label{subsection_1}
\begin{definition}
\label{def_wf_dual}
The $\G$-wave front set and the $\Ginf$-wave front set of a functional $T$ in $\LL(\Gc(\Om),\wt{\C})$ are defined as follows:
\beq
\label{WFGT}
\WF_\G(T):=\hskip-5pt \bigcap_{\substack{a(x,D)\in\,\Oprop{0}(\Om)\\[0.1cm] a(x,D)T\, \in\, \G(\Om)}}\hskip-5pt \compl{\Ellsc(a)},
\eeq
\beq
\label{WFGinfT}
\WF_{\Ginf}(T):=\hskip-5pt \bigcap_{\substack{a(x,D)\in\,\Oprop{0}(\Om)\\[0.1cm] a(x,D)T\, \in\, \Ginf(\Om)}}\hskip-5pt \compl{\Ellsc(a)}.
\eeq
\end{definition}
\begin{remark}
\label{rem_3}
\leavevmode
\begin{trivlist}
\item[(i)] As observed in \cite{GH:05} the action of a pseudodifferential operator with generalized symbol do not change by adding negligible nets of symbols of order $-\infty$. This means that considering the set $\Syscu^{0/-\infty}(\Om\times\R^n) := \Sscu^0(\Om\times\R^n) / \Nuinf(\Om\times\R^n)$ of slow scale generalized symbols of refined order and the corresponding set $\Oprop{0/-\infty}(\Om)$ of properly supported pseudodifferential operators, $\WF_\G(T)$ and $\WF_{\Ginf}(T)$ can be defined equivalently by replacing the set $\Oprop{0}(\Om)$ with $\Oprop{0/-\infty}(\Om)$ in \eqref{WFGT} and \eqref{WFGinfT}. Clearly all the results of Section \ref{section_pseudo} are valid for pseudodifferential operators with generalized symbols of refined order.
\item[(ii)] By standard procedure of lifting symbol orders with $(1-\Delta)^{m/2}$ we easily show that we may take the intersections over $a(x,D)\in\Oprop{m}(\Om)$ (or $\Oprop{m/-\infty}(\Om)$) in both \eqref{WFGT} and \eqref{WFGinfT}.
\item[(iii)] By Theorem 3.10 in \cite{GH:05} it is clear that the notion of $\Ginf$-wave front set coincides with the usual generalized wave front set (see \cite[(3.10)]{GH:05}) on the Colombeau algebra $\G(\Om)$.
\end{trivlist}
\end{remark}
\begin{proposition}
\label{prop_sing_supp}
Let $\pi:\CO{\Om}\to\Om:(x,\xi)\to x$. For any basic functional $T$ in $\LL(\Gc(\Om),\wt{\C})$, 
\beq
\label{pi_1}
\pi(\WF_\G(T)) = \singsupp_\G\,T
\eeq
and
\beq
\label{pi_2}
\pi(\WF_{\Ginf}(T)) = \singsupp_{\Ginf}\,T.
\eeq
\end{proposition}
\begin{proof}
The proof is a revised version of the proof of Proposition 2.8 in \cite{GH:05} by employing the new concepts introduced in the dual context. Crucial are the mapping properties of the generalized pseudodifferential operators acting on $\LL(\Gc(\Om),\wt{\C})$ here involved and the definitions of $\G$- and $\Ginf$-singular support. 

We begin with \eqref{pi_1}, by proving that $\Om\setminus\singsupp_\G\,T\subseteq\Om\setminus\pi(\WF_\G(T))$. If $x_0\in\Om\setminus\singsupp_\G\,T$ then there exists $\phi\in\Cinfc(\Om)$ with $\phi(x_0)=1$ such that $\phi(x,D)T=\phi T\in\G(\Om)$. The multiplication operator $\phi(x,D)$ belongs to $\Oprop{0}(\Om)$ and its symbol is (slow scale) micro-elliptic at $(x_0,\xi_0)$ for all $\xi_0\neq 0$. Therefore, $x_0\not\in\Om\setminus\pi(\WF_\G(T))$.\\
To show the opposite inclusion let $x_0\in\Om\setminus\pi(\WF_\G(T))$. For all $\xi\neq 0$ there exists $a\in\Syscu^{0/-\infty}(\Om\times\R^n)$ slow scale micro-elliptic at $(x_0,\xi)$ such that $a(x,D)$ is properly supported and $a(x,D)T\in\G(\Om)$. Arguing as in the proof of \cite[Proposition 2.8]{GH:05} we find a finite number of generalized symbols $a_i\in\Syscu^{0/-\infty}(\Om\times\R^n)$ such that $a_i(x,D)$ is properly supported and $a_i(x,D)T\in\G(\Om)$. Let $A:=\sum_{i=1}^N a_i(x,D)^\ast a_i(x,D)$. Since $a_i(x,D)T\in\G(\Om)$ and each $a_i(x,D)^\ast$ maps $\G(\Om)$ into $\G(\Om)$ we conclude that $AT\in\G(\Om)$. The arguing at the level of generalized symbols developed for Proposition 2.8 in \cite{GH:05} shows that there exists a slow scale elliptic symbol $b(x,\xi)\in\Syscu^{0/-\infty}(\Om\times\R^n)$ such that $b(x,D)\in\Oprop{0/-\infty}(\Om)$ and $b(x,D)T|_V=AT|_V\in\G(\Om)$ on some neighborhood $V$ of $x_0$. Since $T$ is a basic functional an application of \eqref{pseudoloc_A} in Theorem \ref{theo_reg_dual} leads to $\singsupp_\G\, T=\singsupp_\G\, b(x,D)T$ and consequently $V\cap\singsupp_\G\,T=\emptyset$. This shows that $x_0\not\in\singsupp_\G\,T$.

The proof of the second assertion is immediate. In proving $\Om\setminus\pi(\WF_{\Ginf}(T))\subseteq\Om\setminus\singsupp_{\Ginf}\,T$ is essential to note that since the slow scale generalized symbols can be seen as a special kind of regular symbols, the operators of $\Oprop{0/-\infty}(\Om)$ map $\Ginf(\Om)$ into itself. Hence, if $a(x,D)T\in\Ginf(\Om)$ then $AT\in\Ginf(\Om)$ and $b(x,D)T|_V\in\Ginf(V)$. An application of \eqref{pseudoloc_B} in Theorem \ref{theo_reg_dual} allows to conclude that $\singsupp_{\Ginf}\,T=\singsupp_{\Ginf}\,b(x,D)T$ and completes the proof.
\end{proof}

\subsection{Fourier transform-characterization of $\WF_\G(T)$ and $\WF_{\Ginf}(T)$ when $T$ is a basic functional}
\label{subsection_2}
The Fourier transform-characterization of the generalized wave front sets introduced before needs some preliminary microlocal results concerning the action of a generalized pseudodifferential operator on the dual $\LL(\Gc(\Om),\wt{\C})$. Whereas for regular symbols it is relevant to talk of $\Ginf$-regularity in a conical neighborhood (see \cite[Definition 3.1]{GH:05}), in the larger class of generalized symbols it has a meaning to talk of microlocal $\G$-regularity.
\begin{definition}
\label{def_micro_supp}
Let $a\in\wt{\mathcal{S}}^m_{\rho,\delta}(\Om\times\R^n)$ and $(x_0,\xi_0)\in\CO{\Om}$. The symbol $a$ is $\G$-smoothing at $(x_0,\xi_0)$ if there exist a representative $(a_\eps)_\eps$ of $a$, a relatively compact open neighborhood $U$ of $x_0$ and a conic neighborhood $\Gamma$ of $\xi_0$ such that
\begin{multline}
\label{est_micro_G}
\forall m\in\R\, \forall\alpha,\beta\in\N^n\, \exists N\in\N\, \exists c>0\, \exists\eta\in(0,1]\, \forall(x,\xi)\in U\times\Gamma\, \forall\eps\in(0,\eta]\\
|\partial^\alpha_\xi\partial^\beta_x a_\eps(x,\xi)|\le c\lara{\xi}^m\eps^{-N}.
\end{multline}
The symbol $a$ is $\Ginf$-smoothing at $(x_0,\xi_0)$ if there exist a representative $(a_\eps)_\eps$ of $a$, a relatively compact open neighborhood $U$ of $x_0$, a conic neighborhood $\Gamma$ of $\xi_0$ and a natural number $N\in\N$ such that  
\begin{multline}
\label{est_micro_Ginf}
\forall m\in\R\, \forall\alpha,\beta\in\N^n\, \exists c>0\, \exists\eta\in(0,1]\, \forall(x,\xi)\in U\times\Gamma\, \forall\eps\in(0,\eta]\\
|\partial^\alpha_\xi\partial^\beta_x a_\eps(x,\xi)|\le c\lara{\xi}^m\eps^{-N}.
\end{multline}
We define the \emph{$\G$-microsupport} of $a$, denoted by $\mu\,\supp_\G(a)$, as the complement of the set of points $(x_0,\xi_0)$ where $a$ is $\G$-smoothing and the \emph{$\Ginf$-microsupport} of $a$, denoted by $\mu\,\supp_{\Ginf}(a)$, as the complement of the set of points $(x_0,\xi_0)$ where $a$ is $\Ginf$-smoothing.
\end{definition}
In analogy with \cite{GH:05} when $a\in\wt{\mathcal{S}}^{\,m/-\infty}_{\rho,\delta}(\Om\times\R^n)$ we denote the complements of the sets of points $(x_0,\xi_0)\in\CO{\Om}$ where \eqref{est_micro_G} and \eqref{est_micro_Ginf} hold for some representative of $a$ by $\mu_\G(a)$ and $\mu_{\Ginf}(a)$ respectively. It is clear that:
\begin{itemize}
\item[(i)] if $a\in{\wt{\mathcal{S}}}^{-\infty}(\Om\times\R^n)$ then $\mu\,\supp_\G(a)=\emptyset$;
\item[(ii)] if $a\in{\wt{\mathcal{S}}}^{-\infty}_{\rm{rg}}(\Om\times\R^n)$ then $\mu\,\supp_{\Ginf}(a)=\emptyset$;
\item[(iii)] if $a\in{\wt{\mathcal{S}}}^{\,m/-\infty}_{\rho,\delta}(\Om\times\R^n)$ and $\mu_\G(a)=\emptyset$ then $a\in\wt{\mathcal{S}}^{-\infty}(\Om\times\R^n)$;
\item[(iv)] if $a\in{\wt{\mathcal{S}}}^{\,m/-\infty}_{\rho,\delta,{\rm{rg}}}(\Om\times\R^n)$ and $\mu_{\Ginf}(a)=\emptyset$ then $a\in\wt{\mathcal{S}}^{-\infty}_{\rm{rg}}(\Om\times\R^n)$;
\item[(v)] when $a$ is a classical symbol then $\mu\,\supp(a)=\mu_\G(a)=\mu_{\Ginf}(a)$.
\end{itemize}
In the sequel we work under the hypothesis $0\le\delta<\rho\le 1$. We recall that when $a(x,D)$ and $b(x,D)$ are properly supported pseudodifferential operators with symbols $a\in\wt{\mathcal{S}}^m_{\rho,\delta}(\Om\times\R^n)$ and $b\in\wt{\mathcal{S}}^{m'}_{\rho,\delta}(\Om\times\R^n)$ respectively, then $a(x,D)\circ b(x,D)$ is properly supported itself and has generalized symbol $a\sharp b$ in $\wt{\mathcal{S}}^{m+m'}_{\rho,\delta}(\Om\times\R^n)$. $a\sharp b$ has asymptotic expansion $\sum_\gamma\partial^\gamma_\xi a D^\gamma_x b/\gamma !$ in the sense that for all representatives $(a_\eps)_\eps$ and $(b_\eps)_\eps$ of $a$ and $b$ respectively there exists a representative $((a\sharp b)_\eps)_\eps$ of $a\sharp b$ such that for all $r\in\N\setminus 0$,
\beq
\label{asymp_gen}
\biggl((a\sharp b)_\eps-\sum_{|\gamma|=0}^{r-1}\frac{1}{\gamma !}\partial^\gamma_\xi a_\eps D^\gamma_x b_\eps  \biggr)_\eps\in\mathcal{M}_{S^{m+m'-(\rho-\delta)r}_{\rho,\delta}(\Om\times\R^n)}.
\eeq
Note that when $a$ and $b$ are regular symbols then $a\sharp b$ is regular. More precisely for $K\Subset\Om$, the assumption $|a_\eps|^{(m)}_{\rho,\delta,K,\alpha,\beta}=O(\eps^{-N})$ and $|b_\eps|^{(m)}_{\rho,\delta,K,\alpha,\beta}=O(\eps^{-N'})$, valid for all $\alpha,\beta\in\N^n$, implies that 
\beq
\label{asymp_gen_2}
\biggl|(a\sharp b)_\eps-\sum_{|\gamma|=0}^{r-1}\frac{1}{\gamma !}\partial^\gamma_\xi a_\eps D^\gamma_x b_\eps \biggr|^{(m+m'-(\rho-\delta)r)}_{\rho,\delta,K,\alpha,\beta}=O(\eps^{-N-N'})
\eeq
for all $r\in\N\setminus 0$ and $\alpha,\beta\in\N^n$.
\begin{proposition}
\label{prop_micro_supp_product}
Let $a(x,D)$ and $b(x,D)$ be properly supported pseudodifferential operators with generalized symbols.
\begin{itemize}
\item[(i)] If $a\in\wt{\mathcal{S}}^m_{\rho,\delta}(\Om\times\R^n)$ and $b\in\wt{\mathcal{S}}^{m'}_{\rho,\delta}(\Om\times\R^n)$ then 
\beq
\label{musupp_1}
\mu\, \supp_\G(a\sharp b)\,\subseteq\,\mu\, \supp_\G(a)\,\cap\,\mu\, \supp_\G(b).
\eeq
\item[(ii)] If $a\in\wt{\mathcal{S}}^m_{\rho,\delta,{\rm{rg}}}(\Om\times\R^n)$ and $b\in\wt{\mathcal{S}}^{m'}_{\rho,\delta,{\rm{rg}}}(\Om\times\R^n)$ then 
\beq
\label{musupp_2}
\mu\, \supp_{\Ginf}(a\sharp b)\,\subseteq\,\mu\, \supp_{\Ginf}(a)\,\cap\,\mu\, \supp_{\Ginf}(b).
\eeq
\end{itemize}
\end{proposition}
When we deal with symbols of refined order we have that $\mu\, \supp_\G$ and $\mu\, \supp_{\Ginf}$ can be replaced by $\mu_\G$ and $\mu_{\Ginf}$ respectively in \eqref{musupp_1} and \eqref{musupp_2}.
\begin{proof}
$(i)$ Assume that $(x_0,\xi_0)\not\in\mu\, \supp_\G(a)$. This means that \eqref{est_micro_G} holds for some representative $(a_\eps)_\eps$ of $a$ in a region $U\times\Gamma$. Combined with the properties of $b$  we have that:
\begin{multline}
\label{ab_1}
\forall l\in\R\, \forall\alpha,\beta,\gamma\in\N^n\, \exists N\in\N\, \exists c>0\, \exists\eta\in(0,1]\, \forall (x,\xi)\in U\times\Gamma\, \forall\eps\in(0,\eta]\\
|\partial^\alpha_\xi\partial^\beta_x(\partial^\gamma_\xi a_\eps D^\gamma_x b_\eps)(x,\xi)|\le c\lara{\xi}^l\eps^{-N}.
\end{multline}
By \eqref{asymp_gen}, taking for each $\alpha,\beta\in\N^n$ the integer $r$ large enough such that $m+m'-(\rho-\delta)r-\rho|\alpha|+\delta|\beta|\le l$, we conclude that the following assertion
\begin{multline}
\label{ab_2}
\forall l\in\R\, \forall\alpha,\beta,\gamma\in\N^n\, \exists N\in\N\, \exists c>0\, \exists\eta\in(0,1]\, \forall (x,\xi)\in U\times\Gamma\, \forall\eps\in(0,\eta]\\
|\partial^\alpha_\xi\partial^\beta_x (a\sharp b)_\eps(x,\xi)|\le c\lara{\xi}^l\eps^{-N}
\end{multline}
holds for some representative $((a\sharp b)_\eps)_\eps$ of $a\sharp b$. Hence, $(x_0,\xi_0)\not\in\mu\, \supp_\G(a\sharp b)$.

$(ii)$ When we deal with regular symbols and their $\Ginf$-microsupports, \eqref{ab_1} is transformed in
\begin{multline*}
\exists N,N'\in\N\, \forall l\in\R\, \forall\alpha,\beta,\gamma\in\N^n\, \exists c>0\, 
 \exists\eta\in(0,1]\, \forall (x,\xi)\in U\times\Gamma\, \forall\eps\in(0,\eta]\\
|\partial^\alpha_\xi\partial^\beta_x(\partial^\gamma_\xi a_\eps D^\gamma_x b_\eps)(x,\xi)|\le c\lara{\xi}^l\eps^{-N-N'}.
\end{multline*} 
By means of \eqref{asymp_gen_2} it is immediate to obtain that for all order of derivatives and for all $l\in\R$
\[
|\partial^\alpha_\xi\partial^\beta_x(a\sharp b)_\eps(x,\xi)|\le c\lara{\xi}^l\eps^{-N-N'}
\]
on $U\times\Gamma$ when $\eps$ is small enough, i.e., $(x_0,\xi_0)\not\in\mu\, \supp_{\Ginf}(a\sharp b)$.
\end{proof}
\begin{theorem}
\label{theo_WF_a(x,D)}
Let $a(x,D)$ be a properly supported pseudodifferential operator with generalized symbol and $T$ be a basic functional in $\LL(\Gc(\Om),\wt{\C})$.
\begin{itemize}
\item[(i)] If $a\in\wt{\mathcal{S}}^{m/-\infty}_{\rho,\delta}(\Om\times\R^n)$ then 
\beq
\label{WF(a(x,D))_1}
\WF_\G(a(x,D)T)\,\subseteq\,\WF_\G(T)\,\cap\,\mu_\G(a).
\eeq
\item[(ii)] If $a\in\wt{\mathcal{S}}^{m/-\infty}_{\rho,\delta,{\rm{rg}}}(\Om\times\R^n)$ then 
\beq
\label{WF(a(x,D))_2}
\WF_{\Ginf}(a(x,D)T)\,\subseteq\,\WF_{\Ginf}(T)\,\cap\,\mu_{\Ginf}(a).
\eeq
\end{itemize}
\end{theorem} 
As for Proposition \ref{prop_sing_supp} the proof is an elaboration of the proof of Theorem 3.6 in \cite{GH:05}.
\begin{proof}
$(i)$ We prove the first assertion in two steps.

\emph{Step 1:} \enspace\enspace ${\WF_\G(a(x,D)T)}\subseteq \mu_\G(a)$.

If $(x_0,\xi_0)\not\in\mu_\G(a)$ then \eqref{est_micro_G} holds on some $U\times\Gamma$, and by Lemma 3.4 in \cite{GH:05} we find $q\in{S}^{0}(\Om\times\R^n)\subseteq\Syru^{0/-\infty}(\Om\times\R^n)$, which is micro-elliptic at $(x_0,\xi_0)$ with $\mu\, \supp(q)=\mu_\G(q)\subseteq U\times\Gamma$. Applying Proposition \ref{prop_micro_supp_product} we obtain that $q(x,D)a(x,D)$ is a properly supported pseudodifferential operator with symbol $q\sharp a\in \wt{\mathcal{S}}^{m/-\infty}_{\rho,\delta}(\Om\times\R^n)$ and $\mu_\G(q\sharp a)\subseteq \mu_\G(q)\cap \mu_\G(a)\subseteq (U\times\Gamma) \cap \mu_\G(a) = \emptyset$. This shows that $q\sharp a\in\wt{\mathcal{S}}^{-\infty}(\Om\times\R^n)$ and that $q(x,D)a(x,D)$ has kernel in $\G(\Om\times\Om)$. By Proposition \ref{prop_reg_func}$(ii)$ we have that $q(x,D)a(x,D)T\in\G(\Om)$. Hence, $(x_0,\xi_0)\not\in\WF_\G(a(x,D)T)$.

\emph{Step 2:}\enspace\enspace ${\WF_\G(a(x,D)T)}\subseteq \WF_\G(T)$.

Let $(x_0,\xi_0)\not\in{\WF_\G(T)}$. By definition of $\G$-wave front set there exists $p(x,D)\in\Oprop{0/-\infty}(\Om)$ such that $p$ is slow scale micro-elliptic at $(x_0,\xi_0)$ and $p(x,D)T\in\G(\Om)$. As shown in the proof of \cite[Theorem 3.6]{GH:05} there exist:
\begin{itemize}
\item[-] $r(x,D)\in\Oprop{0/-\infty}(\Om)$ whose symbol $r$ is classical, micro-elliptic at $(x_0,\xi_0)$ and with $\mu\, \supp(r)$ contained in a conic neighborhood $U'\times\Ga'$ of $(x_0,\xi_0)$;
\item[-] $p(x,D)^\ast p(x,D)=\sigma(x,D)\in\Oprop{0/-\infty}(\Om)$;
\item[-] $b(x,D)\in\Oprop{0/-\infty}(\Om)$ whose symbol $b$ is slow scale elliptic and  $\mu_{\Ginf}(b-\sigma)\cap(U'\times\Ga')=\emptyset$;
\item[-] $t(x,D)\in\Opropr{0/-\infty}(\Om)$ parametrix of $b(x,D)$.
\end{itemize}
As a consequence $s(x,D):=r(x,D)a(x,D)t(x,D)$ is a properly supported pseudodifferential operator with symbol in $\wt{\mathcal{S}}^{m/-\infty}_{\rho,\delta}(\Om\times\R^n)$. We write the difference $r(x,D)a(x,D)T-s(x,D)p(x,D)^\ast p(x,D)T$ as 
\beq
\label{diff_12}
r(x,D)a(x,D)\big(T-t(x,D)b(x,D)T\big) + r(x,D)a(x,D)t(x,D)\big(b(x,D)-\sigma(x,D)\big)T.
\eeq
Since $I-t(x,D)b(x,D)$ has kernel in $\Ginf(\Om\times\Om)$ then by Proposition \ref{prop_reg_func}$(iii)$, $T-t(x,D)b(x,D)T\in\Ginf(\Om)$ and by the mapping properties of $r(x,D)$ and $a(x,D)$ the first summand in \eqref{diff_12} belongs to $\G(\Om)$. An iterated application of Proposition \ref{prop_micro_supp_product} stated for symbols of refined order proves that the second summand can be written as the action on $T$ of a properly supported pseudodifferential operator $d(x,D)$ with generalized symbol of refined order $m$ having $\G$-microsupport contained in the region $\mu\, \supp(r)\cap\mu_\G(b-\sigma)\subseteq\mu\, \supp(r)\cap\mu_{\Ginf}(b-\sigma)\subseteq U'\times\Ga'\cap\mu_{\Ginf}(b-\sigma)=\emptyset$. This means that $d\in\wt{\mathcal{S}}^{-\infty}(\Om\times\R^n)$ and by Proposition \ref{prop_charac_infty}$(i)$ and Proposition \ref{prop_reg_func}$(ii)$ we conclude that $d(x,D)T\in\G(\Om)$. Therefore, \eqref{diff_12} gives a generalized function in $\G(\Om)$.

Let us now consider $r(x,D)$. Recalling that $p(x,D)T\in\G(\Om)$ and that the operators $p(x,D)^\ast$ and $s(x,D)$ map $\G(\Om)$ into itself, the considerations above imply that $r(x,D)a(x,D)T\in\G(\Om)$. Thus, $(x_0,\xi_0)\not\in\WF_\G(a(x,D)T)$ since $r$ is micro-elliptic at $(x_0,\xi_0)$.

$(ii)$ When $a\in\wt{\mathcal{S}}^{m/-\infty}_{\rho,\delta,{\rm{rg}}}(\Om\times\R^n)$ then $a(x,D)$ maps $\Ginf(\Om)$ into itself. Since the same mapping property holds for $r(x,D)$ it follows that the first summand in \eqref{diff_12} belongs to $\Ginf(\Om)$. By iterated application of Proposition \ref{prop_micro_supp_product} the second summand $d(x,D)T$ has symbol $d$ with $\Ginf$-microsupport contained in $\mu\, \supp(r)\cap\mu_{\Ginf}(b-\sigma)=\emptyset$, that is $d\in\wt{\mathcal{S}}^{-\infty}_{\rm{rg}}(\Om\times\R^n)$ and the kernel of the corresponding pseudodifferential operator is an element of $\Ginf(\Om\times\Om)$. Proposition \ref{prop_reg_func}$(iii)$ yields $d(x,D)T\in\Ginf(\Om)$ and then \eqref{diff_12} gives a generalized function in $\Ginf(\Om)$. Recalling that by definition of $\Ginf$-wave front set $p(x,D)T\in\Ginf(\Om)$ and that $p(x,D)^\ast$ and $s(x,D)$ map $\Ginf(\Om)$ into itself, we obtain that $r(x,D)(a(x,D)T)\in\Ginf(\Om)$ which implies that $(x_0,\xi_0)\not\in\WF_{\Ginf}(a(x,D)T)$.
\end{proof}
Note that if $a\in\wt{\mathcal{S}}^m_{\rho,\delta}(\Om\times\R^n)$ and $\kappa$ is the quotient map from $\M_{S^m_{\rho,\delta}(\Om\times\R^n)}$ onto $\wt{\mathcal{S}}^{m/-\infty}_{\rho,\delta}(\Om\times\R^n)$ then 
\beq
\label{musuppG_a}
 \mu\, \supp_\G(a) = \bigcap_{(a_\eps)_\eps \in a} \mu_\G\big(\,\kappa\big((a_\eps)_\eps\big)\,\big).
\eeq 
Indeed, for every representative $(a_\eps)_\eps$ of the symbol $a$ we have that $\mu\, \supp_\G(a)\subseteq\mu_\G(\kappa((a_\eps)_\eps))$ and if $(x_0,\xi_0)\not\in\mu\, \supp_\G(a)$ then $(x_0,\xi_0)\not\in\mu_\G(\kappa((a_\eps)_\eps))$ for some $(a_\eps)_\eps\in a$. In the same way 
\beq
\label{musuppGinf_a}
 \mu\, \supp_{\Ginf}(a) = \bigcap_{(a_\eps)_\eps \in a} \mu_{\Ginf}\big(\,\kappa\big((a_\eps)_\eps\big)\,\big).
\eeq
We are ready now to prove the following corollary of Theorem \ref{theo_WF_a(x,D)}.
\begin{corollary}
\label{cor_WF_a(x,D)}
For any properly supported pseudodifferential operator $a(x,D)$ with symbol $a\in\wt{\mathcal{S}}^{m}_{\rho,\delta}(\Om\times\R^n)$ and for any $T$ basic functional in $\LL(\Gc(\Om),\wt{\C})$, 
\beq
\label{WF(a(x,D))_3}
\WF_\G(a(x,D)T)\,\subseteq\,\WF_\G(T)\,\cap\,\mu\, \supp_\G(a).
\eeq
Similarly, if $a\in\wt{\mathcal{S}}^{m}_{\rho,\delta,{\rm{rg}}}(\Om\times\R^n)$ then 
\beq
\label{WF(a(x,D))_4}
\WF_{\Ginf}(a(x,D)T)\,\subseteq\,\WF_{\Ginf}(T)\,\cap\,\mu\, \supp_{\Ginf}(a).
\eeq
\end{corollary}
\begin{proof}
For any representative $(a_\eps)_\eps$ of $a$ the generalized symbol $\kappa((a_\eps)_\eps)=(a_\eps)_\eps+\Neg^{-\infty}(\Om\times\R^n)$ satisfies the hypotheses of Theorem \ref{theo_WF_a(x,D)} and the corresponding operator $\kappa((a_\eps)_\eps)(x,D)$ coincides with $a(x,D)$. Hence from Theorem \ref{theo_WF_a(x,D)}$(i)$ we have that
\[
\bigcap_{(a_\eps)_\eps\in a}\WF_\G(\kappa((a_\eps)_\eps)(x,D)T)\,\subseteq\, \WF_\G(T)\, \cap\, \bigcap_{(a_\eps)_\eps\in a}\mu_\G(\kappa((a_\eps)_\eps)).
\]
Clearly the properties of $\kappa((a_\eps))(x,D)$ and equality \eqref{musuppG_a} lead to 
\[
\WF_\G(a(x,D)T)\,\subseteq\, \WF_\G(T)\, \cap\, \mu\, \supp_\G(a).
\]
The proof of \eqref{WF(a(x,D))_4} when $a$ is a regular symbol is an analogous combination of Theorem \ref{theo_WF_a(x,D)} with \eqref{musuppGinf_a}.
\end{proof}
\begin{remark}
\label{rem_example}
Note that the assumption of regularity of the symbol $a$ is essential in order to get \eqref{WF(a(x,D))_4}. Indeed, for $\phi\in\Cinfc(\R)$ let us consider $\phi(x/\eps)=\eps^{-1}\phi(x/\eps)$, the generalized function $[(\phi_\eps)_\eps]\in\G(\R)$ and the multiplication operator $a(x,D):T\to[(\phi_\eps)_\eps]T$. The generalized symbol determined by $[(\phi_\eps)_\eps]$ is not regular since $[(\phi_\eps)_\eps]\in\G(\R)\setminus\Ginf(\R)$. Taking now the basic functional $T(u)=\int_\R u(x)\, dx$ of $\LL(\Gc(\R),\wt{\C})$ we have that $\WF_{\Ginf}(T)=\emptyset$ while $\WF_{\Ginf}(a(x,D)T)=\WF_{\Ginf}([(\phi_\eps)_\eps])\neq\emptyset$.
\end{remark}
As in the classical theory \cite{Folland:95} we introduce notions of microsupport for operators. In the case of generalized psuedodifferential operators, taking into account the non-injectivity when mapping symbols to operators (cf. \cite{GGO:03}), we distinguish the corresponding notions for symbols and operators.
\begin{definition}
\label{def_microsupp_oper}
Let $A$ be a properly supported pseudodifferential operator with generalized symbol in $\wt{\mathcal{S}}^m_{\rho,\delta}(\Om\times\R^n)$. We define the \emph{$\G$-microsupport} of $A$ as
\[
\Bigmu\,\supp_\G(A):=\bigcap_{\substack{a\in\wt{\mathcal{S}}^m_{\rho,\delta}(\Om\times\R^n)\\ a(x,D)=A}}\mu\,\supp_\G(a).
\]
Let $A$ be a properly supported pseudodifferential operator with generalized symbol in $\wt{\mathcal{S}}^m_{\rho,\delta,{\rm{rg}}}(\Om\times\R^n)$. We define the \emph{$\Ginf$-microsupport} of $A$ as
\[
\Bigmu\,\supp_{\Ginf}(A):=\bigcap_{\substack{a\in\wt{\mathcal{S}}^m_{\rho,\delta,{\rm{rg}}}(\Om\times\R^n)\\ a(x,D)=A}}\mu\,\supp_{\Ginf}(a).
\]
\end{definition} 
Corollary \ref{cor_WF_a(x,D)} can therefore be stated in the following way: for any properly supported pseudodifferential operator $A$ with generalized symbol and for any basic functional $T\in\LL(\Gc(\Om),\wt{\C})$,
\[
\WF_\G(AT)\, \subseteq\, \WF_\G(T)\, \cap\, \Bigmu\,\supp_{\G}(A).
\]
If $A$ has regular generalized symbol then 
\[
\WF_{\Ginf}(AT)\, \subseteq\, \WF_{\Ginf}(T)\, \cap\, \Bigmu\,\supp_{\Ginf}(A).
\]
\begin{corollary}
\label{cor_hyp}
Let $A$ be a properly supported pseudodifferential operator with generalized hypoelliptic symbol. Then for any basic functional $T\in\LL(\Gc(\Om),\wt{\C})$,
\[
\WF_\G(AT)=\WF_\G(T)
\]
and
\[
\WF_{\Ginf}(AT)=\WF_{\Ginf}(T).
\]
\end{corollary}
\begin{proof}
Since $A=a(x,D)$ where $a\in\underline{\wt{\mathcal{S}}}^m_{\,\rho,\delta,{\rm{rg}}}(\Om\times\R^n)$, Corollary \ref{cor_WF_a(x,D)} implies that $\WF_\G(AT)\subseteq\WF_\G(T)$ and $\WF_{\Ginf}(AT)\subseteq\WF_{\Ginf}(T)$. Let $p(x,D)$ be a parametrix for $A$. The symbol $p$ belongs to $\wt{\underline{\mathcal{S}}}^{-l}_{\rho,\delta,{\rm{rg}}}(\Om\times\R^n)$ and from \eqref{parametrix_eq} we have that $p(x,D)AT=T+RT$, where $R$ is an operator with kernel in $\Ginf(\Om\times\Om)$. Proposition \ref{prop_reg_func}$(iii)$ implies that $RT\in\Ginf(\Om)$. Finally an application of Corollary \ref{cor_WF_a(x,D)} to $p(x,D)$ yields $\WF_\G(T)=\WF_\G(p(x,D)AT)\subseteq\WF_\G(AT)$ and $\WF_{\Ginf}(T)=\WF_{\Ginf}(p(x,D)AT)\subseteq\WF_{\Ginf}(AT)$.
\end{proof} 
Note that combining Corollary \ref{cor_hyp} with Proposition \ref{prop_sing_supp} we obtain the equalities between singular supports claimed by Theorem \ref{theo_reg_dual}. Moreover, the statements of the above theorem and corollaries are valid for operators not necessarily properly supported when we consider basic functionals in $\LL(\G(\Om),\wt{\C})$.

Before proving the Fourier transform-characterization of the wave front sets $\WF_\G(T)$ and $\WF_{\Ginf}(T)$ when $T$ is a basic functional in $\LL(\Gc(\Om),\wt{\C})$ we observe that if $\phi\in\Cinfc(\Om)$ then $\phi T$ is a basic functional in $\LL(\G(\R^n),\wt{\C})$ and by Proposition \ref{prop_Fourier_comp} its Fourier transform $\mF(\phi T)$ belongs to $\Gtinf(\R^n)$. In the sequel the regularity of a tempered generalized function is measured on a conic region $\Gamma\subseteq\R^n\setminus 0$ by means of 
\[
\G_{\S,0}(\Gamma):=\{u\in\Gt(\R^n):\ \exists (u_\eps)_\eps \text{repr. of $u$}\ \forall l\in\R\ \exists N\in\N\qquad \sup_{\xi\in\Gamma}\lara{\xi}^l|u_\eps(\xi)|=O(\eps^{-N})\, \text{as $\eps\to 0$}\},
\]
and
\[
\Ginf_{\S\hskip-2pt,0}(\Gamma):=\{u\in\Gt(\R^n):\ \exists (u_\eps)_\eps \text{repr. of $u$}\ \exists N\in\N\ \forall l\in\R\qquad \sup_{\xi\in\Gamma}\lara{\xi}^l|u_\eps(\xi)|=O(\eps^{-N})\, \text{as $\eps\to 0$}\}.
\]
Note that if $(u_\eps)_\eps$ and $(u'_\eps)_\eps$ are two different representatives of $u$ fulfilling the condition which defines $\G_{\S,0}(\Gamma)$ (or $\Ginf_{\S\hskip-2pt,0}(\Gamma)$) then their difference has the property $\sup_{\xi\in\Gamma}\lara{\xi}^l|(u_\eps-u'_\eps)(\xi)|=O(\eps^{q})$ for all $l\in\R$ and all $q\in\N$.
\begin{theorem}
\label{theo_charac_wave}
Let $T$ be a basic functional in $\LL(\Gc(\Om),\wt{\C})$.
\begin{itemize}
\item[(i)] $(x_0,\xi_0)\not\in\WF_\G T$ if and only if there exists a conic neighborhood $\Gamma$ of $\xi_0$ and a cut-off function $\phi\in\Cinfc(\Om)$ with $\phi(x_0)=1$ such that $$\mF(\phi T)\in\G_{\S,0}(\Gamma).$$
\item[(ii)] $(x_0,\xi_0)\not\in\WF_{\Ginf} T$ if and only if there exists a conic neighborhood $\Gamma$ of $\xi_0$ and a cut-off function $\phi\in\Cinfc(\Om)$ with $\phi(x_0)=1$ such that $$\mF(\phi T)\in\Ginf_{\S\hskip-2pt,0}(\Gamma).$$
\end{itemize}
\end{theorem}
\begin{proof}
$(i)$ We first prove that if $(x_0,\xi_0)$ is a point in $\CO{\Om}$ such that $\mF(\phi T)\in\G_{\S,0}(\Ga)$ for some conic neighborhood $\Ga$ of $\xi_0$ and some cut-off function $\phi$ with $\phi(x_0)=1$, then $(x_0,\xi_0)\not\in\WF_\G(T)$. As noted in \cite[Remark 3.5]{GH:05} there exists $p(\xi)\in{S}^0(\Om\times\R^n)$ with $\supp(p)\subseteq\Gamma$, which is identically $1$ in a conical neighborhood $\Gamma'$ of $\xi_0$ when $|\xi|\ge 1$. Taking a proper cut-off $\chi$, we can write the properly supported pseudodifferential operator with amplitude $\chi(x,y)p(\xi)\phi(y)$ in the form $\sigma(x,D)\in\Oprop{0/-\infty}(\Om)$, where $\sigma(x,\xi)-p(\xi)\phi(x)\in{S}^{-1}(\Om\times\R^n)$; in particular for any $S\in\LL(\Gc(\Om),\wt{\C})$,  $\sigma(x,D)S-p(D)(\phi S)$ can be seen as the action of a pseudodifferential operator with kernel in $\Cinf(\Om\times\Om)$ on the functional $\phi S\in\LL(\G(\Om),\wt{\C})$. Hence, if $S$ is basic, 
\beq
\label{I_eq}
\sigma(x,D)S-p(D)(\phi S)\in\Ginf(\Om).
\eeq
By assumption the symbol $\sigma$ is micro-elliptic at $(x_0,\xi_0)$ and by definition of Fourier transform on $\LL(\GS(\R^n),\wt{\C})$ we have that for all $u\in\Gc(\R^n)$,
\begin{multline*}
p(D)(\phi T)(u)=\phi T({\ }^tp(D)(u))=\phi T(\mF(p\mF^{\ast}u))=\mF(\phi T)(p\mF^{\ast}u)=\\
\int_\Om\int_{\R^n}e^{iy\xi}\mF(\phi T)(\xi)p(\xi)\, \dslash\xi\, u(y)\, dy,
\end{multline*}
where $p\mF^{\ast}u\in\GS(\R^n)$ and from the hypothesis on $\mF(\phi T)$ it follows that the integral \beq
\label{int_WF}
\int_{\R^n}e^{iy\xi}\mF(\phi T)(\xi)p(\xi)\, \dslash\xi
\eeq
defines a generalized function in $\G(\Om)$. As a consequence $p(D)(\phi T)\in\G(\Om)$ and by \eqref{I_eq} we conclude that $\sigma(x,D)T\in\G(\Om)$. Therefore, $(x_0,\xi_0)\not\in\WF_\G(T)$.

Conversely, suppose $(x_0,\xi_0)\not\in\WF_\G(T)$. There is an open neighborhood $U$ of $x_0$ such that $(x,\xi_0)\in\compl{\WF_\G(T)}$ for all $x\in U$. Choose $\phi\in\Cinf_{\rm{c}}(U)$ with $\phi(x_0)=1$ and define 
\[
\Sigma_\G :=\{ \xi\in\R^n\setminus 0:\ \ \exists x\in\Om\ (x,\xi)\in\WF_\G(\phi T)\}.
\]
By Corollary \ref{cor_WF_a(x,D)} we have that $\WF_\G(\phi T)\subseteq \WF_\G(T) \cap (\supp(\phi)\times\R^n\setminus 0)$ and therefore $\xi_0\notin\Sigma_\G$. Arguing as in the proof of Theorem 3.10 in \cite{GH:05} we find $p(\xi)\in{S}^0(\Om\times\R^n)$ such that $0\le p\le 1$, $p(\xi)=1$ in a conic neighborhood $\Gamma$ of $\xi_0$ when $|\xi|\ge 1$ and $p(\xi)=0$ in a conic neighborhood $\Sigma_0$ of $\Sigma_\G$. By construction $\mu\,\supp(p) \cap (\Om \times {\Sigma_0}) = \emptyset$ and $\WF_\G(\phi T)\subseteq\Omega\times\Sigma_\G$. Therefore, $\WF_\G(p(D)\phi T)\subseteq\WF_\G(\phi T)\cap \mu\,\supp(p)=\emptyset$ and by Proposition \ref{prop_sing_supp}, $p(D)(\phi T)\in\G(\Om)$. Note that the pseudodifferential operator $p(D)$ maps $\GS(\R^n)$ into itself and can be extended to the dual $\LL(\GS(\R^n),\wt{\C})$. Since $\phi T\in\LL(\G(\Om),\wt{\C})\subseteq\LL(\G(\R^n),\wt{\C})\subseteq\LL(\GS(\R^n),\wt{\C})$, $p(D)(\phi T)$ can be also viewed as a basic functional in $\LL(\GS(\R^n),\wt{\C})$ which restricted to $\Gc(\Om)$ is a generalized function of $\G(\Om)$. We now want to study the action of $p(D)(\phi T)$ on a generalized function $u\in\GS(\R^n)$. First of all for $\delta>0$ we define 
\[
\overline{B}_\delta:=\{x\in\R^n:\, {\rm{dist}}(x,\supp\, \phi)\le \delta\}.
\]
Recalling that $\check{p}$ is a Schwartz function outside the origin, i.e., $\sup_{|x|>\lambda}|x^\alpha\partial^\beta\check{p}(x)|<\infty$ for all $\lambda>0$ and $\alpha,\beta\in\N^n$  (\cite[Theorem (8.8a)]{Folland:95}), from the properties of a net $(T_\eps)_\eps\in\E'(\Om)^{(0,1]}$ defining $T$ we have that there exist $K\Subset\Om$, $j\in\N$ and $N\in\N$ such that for all $l\in\R$ and $\alpha\in\N^n$,
\begin{multline}
\label{est_T_p}
\sup_{x\in\R^n\setminus\overline{B}_\delta}|\lara{x}^l\partial^\alpha_x(\phi T_\eps\ast\check{p})(x)|=\sup_{x\in\R^n\setminus\overline{B}_\delta}|T_\eps(y\to\lara{x}^l\partial^\alpha_x\check{p}(x-y)\phi(y))|\\
\le c\eps^{-N}\sup_{x\in\R^n\setminus\overline{B}_\delta}\sup_{y\in K, |\beta|\le j}|\lara{x}^l\partial^\beta_y(\partial^\alpha_x\check{p}(x-y)\phi(y))|\le c'\eps^{-N}.
\end{multline}
when $\eps$ is small enough. Since $N$ depends only on $(T_\eps)_\eps$ we conclude that $(\phi T_\eps\ast \check{p})_\eps\in\ESinf(\R^n\setminus\overline{B}_\delta)$. Take now $\delta$ small enough such that $\overline{B}_{2\delta}$ is a compact set contained in $\Om$ and a covering $(\Om_j)_{j\in\N}$ of $\R^n$ such that $\Om_0=B_{2\delta}$ and $\Om_j\subseteq\R^n\setminus{\overline{B}_\delta}$ for all $j\ge 1$. Let $(\varphi_j)_j$ be a partition of unity subordinated to the covering $(\Om_j)_j$ ($\supp\, \varphi_j\subseteq\Om_j$) fulfilling the following condition: 
\beq
\label{cond_part}
\forall\alpha\in\N^n\, \exists A_\alpha>0\, \forall j\in\N\qquad \sup_{x\in\R^n}|\partial^\alpha\varphi_j(x)|\le A_\alpha 2^{-j|\alpha|}
\eeq
(see \cite[Theorem 6.1]{Wong:99} for details). Making use of this technical tool we complete the proof of the first assertion of the theorem showing that the Fourier transform of $p(D)(\phi T)\in\LL(\GS(\R^n),\wt{\C})$ belongs to $\GS(\R^n)$. From \eqref{cond_part} we have that if $f\in\S(\R^n)$ then $\sum_j\varphi_j f$ converges to $f$ in $\S(\R^n)$. As a consequence for all $u\in\GS(\R^n)$ we can write
\[
p(D)(\phi T)(u)=p(D)(\phi T)(\varphi_0 u)+\biggl[\biggl(\sum_{j=1}^\infty\int_{\R^n}(\phi T_\eps\ast\check{p})(x)\varphi_j(x)u_\eps(x)\, dx\biggr)_\eps\biggr].
\]
Since $\supp\, \varphi_0\subseteq B_{2\delta}\subseteq\overline{B}_{2\delta}\Subset\Om$ the generalized function $\varphi_0 u$ belongs to $\Gc(\Om)$ and by the properties of $p(D)(\phi T)$ discussed above there exists $v\in\G(\Om)$ such that $p(D)(\phi T)(\varphi_0 u)=\int_\Om v(x)\varphi_0u(x)\, dx$. A combination of the estimate \eqref{est_T_p} with the fact that $\supp\, \varphi_j\subseteq\Om_j\subseteq\R^n\setminus \overline{B}_\delta$ for all $j\ge 1$ and the convergence property \eqref{cond_part} allows to conclude that 
\[
\sum_{j=1}^\infty\int_{\R^n}(\phi T_\eps\ast\check{p})(x)\varphi_j(x)u_\eps(x)\, dx =\int_{\R^n}\sum_{j=1}^\infty(\phi T_\eps\ast\check{p})(x)\varphi_j(x)u_\eps(x)\, dx= \int_{\R^n}w_{2,\eps}(x)u_\eps(x)\, dx,
\]
where $(w_{2,\eps})_\eps\in\ESinf(\R^n)$. In other words there exists $w\in\GS(\R^n)$ of the form $w=w_1+w_2$ with $w_1\in\Gc(\R^n)$ and $w_2\in\GSinf(\R^n)$ such that 
\beq
\label{eq_pD}
p(D)(\phi T)(u)=\int_{\R^n}w(x)u(x)\, dx
\eeq
for all $u\in\GS(\R^n)$. At this point by applying the Fourier transform on both the members of \eqref{eq_pD} we arrive at
\[
\int_{\R^n}\mF w(\xi)u(\xi)\, d\xi=\mF(p(D)(\phi T))(u)=\phi T({\, }^tp(D)\mF u)=\phi T(\mF(pu))=\int_{\R^n}\mF(\phi T)(\xi)p(\xi)u(\xi)\, d\xi,
\]
where by Proposition \ref{prop_Fourier_comp} $\mF(\phi T)\in\Gtinf(\R^n)$ and the generalized functions $\mF w$ and $pu$ are elements of $\GS(\R^n)$. Therefore, $\mF(\phi T)p\in\GS(\R^n)$ and by construction of the symbol $p$ it is clear that $\mF(\phi T)\in\G_{\S,0}(\Ga)$.

$(ii)$ The sufficiency of the second assertion is proven as in the case of the first assertion by simply observing that when $(x_0,\xi_0)$ is a point in $\CO{\Om}$ such that $\mF(\phi T)\in\Ginf_{\S\hskip-2pt,0}(\Gamma)$ then the integral in \eqref{int_WF} defines a generalized function in $\Ginf(\Om)$. Analogously, when $(x_0,\xi_0)\not\in\WF_{\Ginf}(T)$ then $p(D)(\phi T)\in\LL(\GS(\R^n),\wt{\C})$ acts on $\Gc(\Om)$ as a generalized function in $\Ginf(\Om)$. It follows that for all $u\in\GS(\R^n)$, $p(D)(\phi T)(u)=\int_{\R^n}w(x)u(x)\, dx$, where $w=w_1+w_2\in\GSinf(\R^n)$ since $w_1\in\Gcinf(\R^n)$ and $w_2\in\GSinf(\R^n)$. Hence, $\mF(\phi T)\in\Ginf_{\S\hskip-2pt,0}(\Gamma)$.
\end{proof}
As for the generalized wave front set of a Colombeau function (cf. \cite[Theorem 3.12]{GH:05}), the $\G$-wave front set and the $\Ginf$-wave front set of a basic functional in $\LL(\Gc(\Om),\wt{\C})$ can be defined by considering only classical pseudodifferential operators in \eqref{WFGT} and \eqref{WFGinfT}. This is already partially proved in the proof of Theorem \ref{theo_charac_wave}. In the sequel we set
\beq
\label{WGcl}
{\rm{W}}_{{\rm{cl}},\G}(T):=\bigcap_{AT\in\G(\Om)}\Char(A)
\eeq
and
\beq
\label{WGinfcl}
{\rm{W}}_{{\rm{cl}},\Ginf}(T):=\bigcap_{AT\in\Ginf(\Om)}\Char(A)
\eeq
where the intersections are taken over all the classical properly supported operators $A\in\Psi^0(\Om)$ such that $AT\in\G(\Om)$ in \eqref{WGcl} and $AT\in\Ginf(\Om)$ in \eqref{WGinfcl}.
\begin{proposition}
\label{pro_cla}
For all basic functionals $T\in\LL(\Gc(\Om),\wt{\C})$,
\[
{\rm{W}}_{{\rm{cl}},\G}(T)=\WF_\G(T)
\]
and 
\[
{\rm{W}}_{{\rm{cl}},\Ginf}(T)=\WF_{\Ginf}(T).
\] 
\end{proposition}
\begin{proof}
The inclusions $\WF_\G(T)\subseteq{\rm{W}}_{{\rm{cl}},\G}(T)$ and $\WF_{\Ginf}(T)\subseteq{\rm{W}}_{{\rm{cl}},\Ginf}(T)$ are obvious. Let now $(x_0,\xi_0)$ be a point in the complement of $\WF_\G(T)$. As in the proof of Theorem \ref{theo_charac_wave} one can find a properly supported operator $P\in\Psi^0(\Om)$ such that $PT\in\G(\Om)$ and $(x_0,\xi_0)\not\in\Char P$. Hence, $(x_0,\xi_0)\not\in{\rm{W}}_{{\rm{cl}},\G}(T)$. In the same way if $(x_0,\xi_0)\not\in\WF_{\Ginf}(T)$ then $PT\in\Ginf(\Om)$ and $(x_0,\xi_0)\not\in{\rm{W}}_{{\rm{cl}},\Ginf}(T)$.
\end{proof}
\begin{remark}
\label{rem_distr}
Let $\iota$ be the embedding of $\D'(\Om)$ into $\G(\Om)$. We denote the composition of $\iota$ with the embedding of $\G(\Om)$ into the dual $\LL(\Gc(\Om),\wt{\C})$ by $\iota'$ and the straightforward embedding $w\to(u\to [(w(u_\eps))_\eps])$ of $\D'(\Om)$ into $\LL(\Gc(\Om),\wt{\C})$ by $\iota_d$. Note that by \cite[Proposition 3.10]{Garetto:05b} for any classical properly supported pseudodifferential operator $A$ and for any distribution $w$ 
the relation $A(\iota_d(w))=\iota_d(Aw)\in\Ginf(\Om)$ implies $Aw\in\Cinf(\Om)$. Hence a combination of this fact with the previous proposition and Remark \ref{rem_3}$(iii)$ yields 
\[
\WF_{\Ginf}(\iota_d(w))={\rm{W}}_{{\rm{cl}},\Ginf}(\iota_d(w))=\WF(w)=\WF_{\Ginf}(\iota(w))=\WF_{\Ginf}(\iota'(w)).
\]

\end{remark}

\section{Noncharacteristic $\G$ and $\Ginf$-regularity}
\label{section_nonch}
The classical result on noncharacteristic regularity for distributional solutions of arbitrary pseudodifferential equations (with smooth symbols) had been extended to generalized pseudodifferential operators with slow scale generalized symbols and Colombeau solutions in \cite[Theorem 4.1]{GH:05}. We conclude the paper by providing a suitable adaptation and extension of this result to the context of basic functionals in $\LL(\Gc(\Om),\wt{\C})$.
\begin{theorem}
\label{theo_non_charac}
If $P=p(x,D)$ is a properly supported pseudodifferential operator with symbol $p\in\Syscu^m(\Om\times\R^n)$ and $T$ is a basic functional in $\LL(\Gc(\Om),\wt{\C})$ then 
\beq
\label{non_1}
\WF_\G(PT)\,\subseteq\, \WF_\G(T)\, \subseteq\, \WF_\G(PT)\, \cup\, \Ellsc(p)^{\rm{c}}
\eeq
and
\beq
\label{non_2}
\WF_{\Ginf}(PT)\,\subseteq\, \WF_{\Ginf}(T)\, \subseteq\, \WF_{\Ginf}(PT)\, \cup\, \Ellsc(p)^{\rm{c}}.
\eeq 
\end{theorem}
\begin{proof}
From Corollary \ref{cor_WF_a(x,D)} the first inclusions in \eqref{non_1} and \eqref{non_2} are clear. Assume now that $(x_0,\xi_0)\not\in\WF_\G(PT)$ and that $p$ is slow scale micro-elliptic there. By definition of $\G$-wave front set we find $a(x,D)\in\Oprop{0}(\Om)$ such that $a(x,D)(p(x,D)T)\in\G(\Om)$. By the (slow scale) symbol calculus and Remark \ref{rem_slow_ell}(ii) we obtain that $a(x,D)p(x,D)$ has a slow scale symbol micro-elliptic at $(x_0,\xi_0)$. Therefore, $(x_0,\xi_0)\not\in\WF_\G(T)$. Analogously, $\WF_{\Ginf}(PT)^{\rm{c}}\,\cap\,\Ellsc(p)\,\subseteq\,\WF_{\Ginf}(T)^{\rm{c}}$.
\end{proof}

\paragraph{Acknowledgement:} The author is grateful to Prof. Michael Oberguggenberger for several
inspiring discussions on the subject.  

\bibliographystyle{abbrv}
\newcommand{\SortNoop}[1]{}

\end{document}